\documentclass[11pt,a4paper]{article}
\usepackage[T1]{fontenc}
\usepackage{lmodern,mathtools,amssymb,amsthm,microtype,graphicx,booktabs}
\usepackage{enumitem}
\usepackage{longtable}
\usepackage{float}
\usepackage[margin=29mm]{geometry}
\usepackage[font=small,labelfont=bf]{caption}
\usepackage{cite,xcolor}
\usepackage{tabularx}
\usepackage[colorlinks=true,linkcolor=blue!60!black,citecolor=blue!60!black,urlcolor=blue!60!black]{hyperref}
\hypersetup{pdftitle={Shape optimization under minimal width and constant width constraints},pdfauthor={Beniamin Bogosel}}
\numberwithin{equation}{section}
\newtheorem{theorem}{Theorem}[section]
\newtheorem{lemma}[theorem]{Lemma}
\newtheorem{proposition}[theorem]{Proposition}
\newtheorem{corollary}[theorem]{Corollary}

\theoremstyle{remark}\newtheorem{remark}[theorem]{Remark}
\newcommand{\R}{\mathbb R}
\newcommand{\dd}{\,\mathrm d}
\newcommand{\E}{\mathcal E}
\newcommand{\J}{\mathcal J}
\newcommand{\ind}{\mathbf1}
\DeclareMathOperator{\conv}{conv}
\DeclareMathOperator{\bh}{bh}
\DeclareMathOperator{\Int}{int}
\title{Shape optimization under minimal width\\and constant width constraints}
\author{Beniamin Bogosel\thanks{Faculty of Exact Sciences, Aurel Vlaicu
University of Arad, 2 Elena Dr\u{a}goi Street, Arad, Romania.
Corresponding author: \href{mailto:beniamin.bogosel@uav.ro}{\texttt{beniamin.bogosel@uav.ro}}.}}
\date{}
\begin{document}
\maketitle
\begin{abstract}
In this paper we study shape optimization problems for planar convex
sets with prescribed minimal width or constant width. We prove that the
Reuleaux triangle minimizes torsional rigidity among sets of constant
width. We also show that the first Dirichlet-Laplace eigenvalue is
maximized by the equilateral triangle under a minimal width constraint
and by the Reuleaux triangle under a constant width constraint.
The proofs use subsets consisting of a disk and three straight or
circular caps. Reflection arguments reduce the study to equilateral
configurations, which contain a polygon or a disk polygon depending on
one parameter. The same cap constructions give a common proof of the
known area, Cheeger and inner parallel area inequalities, and yield
sharp bounds for all positive moments of the distance to the boundary.
We also prove sharp inequalities for concentric disk intersections and
decreasing radial integrals, and compare interaction energies for kernels
close to a positive constant. This includes a nonempty interval of
positive Riesz exponents.
Finally, we investigate whether inclusion monotonicity and Minkowski
concavity suffice to obtain these comparisons for other functionals.
\end{abstract}

\section{Introduction}\label{sec:introduction}
We are interested in minimizing area and torsional rigidity, and
maximizing the Cheeger constant and the first Dirichlet-Laplace
eigenvalue, under two constraints on the width of a planar convex set.
For a planar convex body $K$ with support function
$h_K(e)=\max_{x\in K}x\cdot e$, its minimal width is
\[
 w(K)=\min_{e\in\mathbb S^1}\{h_K(e)+h_K(-e)\}.
\]
It has constant width $w$ if the expression in braces equals $w$ in
every direction. The respective optimizers are the equilateral triangle
and the Reuleaux triangle.
The proofs use a common geometric reduction: each admissible body
contains a disk with three straight or circular caps. Inclusion
monotonicity reduces the functional comparisons to these families,
which are described by a few geometric parameters. For torsion and
the first eigenvalue, reflection arguments and shape derivatives then
determine the optimal cap positions and compare the resulting
equilateral configurations with the reference bodies.

Let $\mathcal T$ be an equilateral triangle of altitude one and
$\mathcal R$ a Reuleaux triangle of constant width one. We use
\[
 -\Delta u_K=1\text{ in }\Int K,\quad u_K\in H^1_0(\Int K),\qquad
 T(K)=\int_Ku_K=\int_K|\nabla u_K|^2,
\]
and write $\lambda_1(K)$ for the first Dirichlet eigenvalue of
$-\Delta$, and $h(K)$ for the Cheeger constant. 
We denote the inradius by $r(K)$ and write
$K_{-t}=\{x:x+\overline B_t\subset K\}$ for the inner parallel set at
distance $t\geq0$, where $\overline B_t$ is the closed disk of radius
$t$ centered at the origin.

\begin{theorem}\label{thm:main}
Let $K$ be a planar convex body and $w>0$. The following sharp
comparisons hold:
\begin{center}\normalfont
\begin{tabular}{@{}lll@{}}
\toprule
Functional & $w(K)=w$ & constant width $w$\\
\midrule
Area & $|K|\geq w^2|\mathcal T|$ & $|K|\geq w^2|\mathcal R|$\\[2pt]
Cheeger constant & $h(K)\leq w^{-1}h(\mathcal T)$ & $h(K)\leq w^{-1}h(\mathcal R)$\\[2pt]
Torsional rigidity & $T(K)\geq w^4T(\mathcal T)$ & $T(K)\geq w^4T(\mathcal R)$\\[2pt]
First eigenvalue & $\lambda_1(K)\leq4\pi^2w^{-2}$ & $\lambda_1(K)\leq w^{-2}\lambda_1(\mathcal R)$\\
\bottomrule
\end{tabular}
\end{center}
Equality in any entry characterizes the corresponding triangle up to
rigid motion and scaling. Moreover, for $\mathcal Q=\mathcal T$ or
$\mathcal R$ in the respective class,
\begin{equation}\label{eq:main-inner}
 |K_{-t}|\geq w^2|\mathcal Q_{-t/w}|\qquad(t\geq0).
\end{equation}
Equality in \eqref{eq:main-inner} at any
$0\leq t\leq wr(\mathcal Q)$ also characterizes $w\mathcal Q$.
\end{theorem}

The new PDE conclusions are the constant-width torsion minimum and
the first-eigenvalue maxima in both width classes. The minimal-width
torsion comparison was proved in \cite{widthpaper}; the area, Cheeger
and inner-parallel-area comparisons are earlier results, recovered
here within the same framework.

Some of the constants in Theorem~\ref{thm:main} are explicit:
\begin{equation}\label{eq:triangle-values}
 |\mathcal T|=\frac1{\sqrt3},\quad
 h(\mathcal T)=3+\sqrt{\pi\sqrt3},\quad
 T(\mathcal T)=\frac{\sqrt3}{180},\quad
 \lambda_1(\mathcal T)=4\pi^2,\qquad
 |\mathcal R|=\frac{\pi-\sqrt3}{2}.
\end{equation}
The first eigenvalue of the equilateral triangle is classical; see \cite{McCartin}.
For torsion, if $\ell_i$ are the three distances to its sides then
$\sum_i\ell_i=1$ and $u=\ell_1\ell_2\ell_3$ solves $-\Delta u=1$.
Integrating the barycentric monomial gives $T=|\mathcal T|/60$.
The Cheeger value follows from the inner-area equation below.

The same subsets also yield integral comparisons. With the origin at an
incenter, we prove
\[
 |K\cap B_R|\geq|\mathcal Q\cap B_R|\qquad(R\geq0)
\]
at width one, with rigidity for $R>r(\mathcal Q)$.
This gives a minimum for every nonnegative decreasing radial weight and
also proves the shape comparison after maximizing the integral over its center.
For two-point energies, angular balancing holds for convex decreasing
kernels. A quantitative estimate within the cap family gives global
minimality for kernels close to a positive constant, including a
nonempty interval of Riesz exponents in both width classes.

\subsection{Previous results}
P\'al's isominwidth theorem \cite{Pal} and the Blaschke--Lebesgue
theorem \cite{Blaschke,Lebesgue} give the area row. The classical
geometry of constant width, including the relations between its
inscribed and circumscribed circles, is treated in Yaglom and
Boltyanskii \cite[Chapter~7]{YaglomBoltyanskii} and in
\cite{ConstantWidthBook}. Other approaches to Blaschke--Lebesgue
include \cite{Chakerian,Harrell,MalagoliControl,Hynd,BogoselMixed,BogoselVariations}.
Ball and spindle convexity provide the curved counterpart of ordinary
convex hulls \cite{BallPolyhedra}.

The straight-cap construction is described in
\cite[Chapter~6]{YaglomBoltyanskii} and used by Lucardesi and Zucco
\cite[Proposition~3.3]{LucardesiZucco} in their quantitative versions of
P\'al's inequality. Their Lemma~3.2 also gives the area monotonicity for
this family. The same subset is used for torsion in \cite{widthpaper}.
The curved three-cap construction and its area comparison already
occur in M.~Bezdek's proof for disk polygons \cite[Section~4]{BezdekDisk}.
More recently Fodor, Robock and Sagmeister \cite{FodorRobockSagmeister}
connected the P\'al and Blaschke--Lebesgue inequalities through ball
convexity, using three caps and a disk subhexagon. We use these
constructions below to study the Dirichlet functionals.
An earlier analytic approach to minimum area under width and curvature
constraints is given by Malagoli \cite{MalagoliWidth}.
Henrot and Lucardesi \cite{HenrotLucardesiAnnulus}
also solved the distinct area problem with both constant width and
inradius prescribed; our intermediate cap bodies need not themselves
satisfy the width constraint.

Henrot and Lucardesi \cite{HenrotLucardesiCheeger} proved the
constant-width Cheeger inequality. The inner-parallel-area and
inner-perimeter comparisons of \cite{BogoselCheeger} give another
geometric proof. Under a minimal-width constraint, Ftouhi, Lucardesi
and Saracco \cite{FtouhiLucardesiSaracco} proved the Cheeger inequality
and its local Hausdorff stability. Their Corollary~4.5 also proves the
inner-parallel-area comparison with its strict equality case.
The straight-cap method in \cite{widthpaper} gives another proof of the
Cheeger inequality and establishes the torsion inequality. Here explicit
erosion of the two cap families gives a common proof of the known
inner-area comparisons. Sharp inequalities for distance moments follow
by integration.
The inner-perimeter result in the constant-width class requires the
additional geometric argument of \cite{BogoselCheeger}.

For nonlocal functionals, Bogosel, Bucur and Fragal\`a
\cite{BogoselBucurFragala} prove a polygonal Hardy--Littlewood inequality
and study symmetry and symmetry breaking for Riesz energies at fixed
area and number of sides. We consider minimum problems under width
constraints, using the same reduction to indicator kernels for one-point
integrals.

For torsion and the first eigenvalue, Philippin \cite{Philippin} proved
the comparisons among regular Reuleaux polygons. Numerical
support-function and convexity-constrained methods
\cite{AntunesBogosel,BogoselNumerical} suggest the same extrema more
generally. The relationship with diameter-constrained spectral
optimization is discussed in \cite{BHL}, where the Reuleaux triangle
was already proposed as a maximizer of $\lambda_1$ at fixed constant
width. Our first-eigenvalue proof
settles the $p=2$ case of the constant-width conjecture formulated in
\cite[Section~5]{HenrotLucardesiCheeger} and the first-eigenvalue case
of the minimal-width conjecture in \cite{widthpaper}.
It uses positivity and log-concavity of the ground state.
Boundary-flux reflection arguments are also related to the triangle
rigidity results of Fragal\`a and Velichkov \cite{FV}. The roles of
inner parallel sets and distance-dependent test functions in torsion
are developed in \cite{FragalaGazzolaLamboley}.

The four functionals, after their natural degree-one normalization,
are Minkowski concave \cite{bm_eigen,BorellTorsion,Colesanti,FtouhiParallel}.
Concavity and indecomposability have a substantial role in shape
optimization \cite{BucurFragalaLamboley}. For constant width, the planar
extreme bodies were characterized by Kallay \cite{Kallay}; related
questions in dimension three are studied in \cite{HyndExtreme}.
It is therefore natural to ask whether these properties imply the
comparisons in Theorem~\ref{thm:main}. We show that this is not the
case, by constructing counterexamples which also satisfy rigid-motion
invariance and strict inclusion monotonicity.

\subsection{Outline of the proof}
Section~\ref{sec:cap-subsets} constructs the cap subsets and their rotation
velocities. Section~\ref{sec:geometric-inequalities} treats inner parallel
sets and Minkowski concavity. Sections~\ref{sec:torsion} and~\ref{sec:eigenvalues}
complete the torsion and first-eigenvalue proofs, respectively.
Integral functionals are treated in Section~\ref{sec:integral-kernels}.
Appendices~\ref{app:variations} and~\ref{app:straight} justify the
first-variation formulas, including admissible motions at cap contacts.

By scaling we may assume that $w=1$. In the following we write
$e_a=(\cos a,\sin a)$ and
$J(x,y)=(-y,x)$, and set $r_0=1-1/\sqrt3$.
At minimal width the original inradius satisfies $1/3\leq r\leq1/2$.
There is a subset $S_r^{\rm s}\subset K$ consisting of $\overline B_r$
and three disjoint straight caps whose vertices have distance $1-r$
from the disk center. At constant width the inradius satisfies
$r_0\leq r\leq1/2$, and the concentric circumcircle has radius $1-r$.
Three contact points with the circumcircle give a subset $S_r^{\rm c}\subset K$ obtained
by adding circular caps of generating radius one. See Figure \ref{fig:proof-overview} for illustrations.

The areas of these cap bodies and of their inner parallel sets are
independent of the angular positions of the three caps. Computing their
derivatives with respect to $r$ gives the area and Cheeger comparisons.
For torsion and the first eigenvalue consider the following constructions:
\begin{itemize}[noitemsep,topsep=0pt]
	\item $E_r$ is the equilateral cap body, with the three cap axes at equal angular distances;
	\item $H_r\subset E_r$ is the alternating hexagon, with vertices on six rays separated by $\pi/3$, at alternating distances $r$ and $1-r$ from the center.
\end{itemize}  
The hexagon has straight edges at minimal width and unit circular edges
at constant width. For minimal width, the disk bounds with cutoffs
\eqref{eq:torsion-cutoff} and \eqref{eq:straight-cutoffs} handle larger
inradii. The angular comparisons are needed only below these cutoffs,
within $1/3<r<\sqrt2-1$; for curved caps they hold throughout
$r_0<r<1/2$. To compare $S_r$ and $E_r$, we use reflections of the
corresponding solutions and move the caps until the angular
configuration is equilateral. The comparison between $H_r$ and
$\mathcal Q$ follows from the shape derivative: we pair reflected
boundary points and compare both their normal velocities and the
normal derivatives of the solutions. The sets $H_r$ are not nested as
$r$ varies, and their inradius may differ from $r$. These constructions
give the comparison chains
\begin{align}
 T(K)&\geq T(S_r)\geq T(E_r)\geq T(H_r)\geq T(\mathcal Q),
 \label{eq:chain}\\
 \lambda_1(K)&\leq\lambda_1(S_r)\leq\lambda_1(E_r)
 \leq\lambda_1(H_r)\leq\lambda_1(\mathcal Q).
 \label{eq:spectral-chain}
\end{align}
\begin{center}
\includegraphics[width=\linewidth]{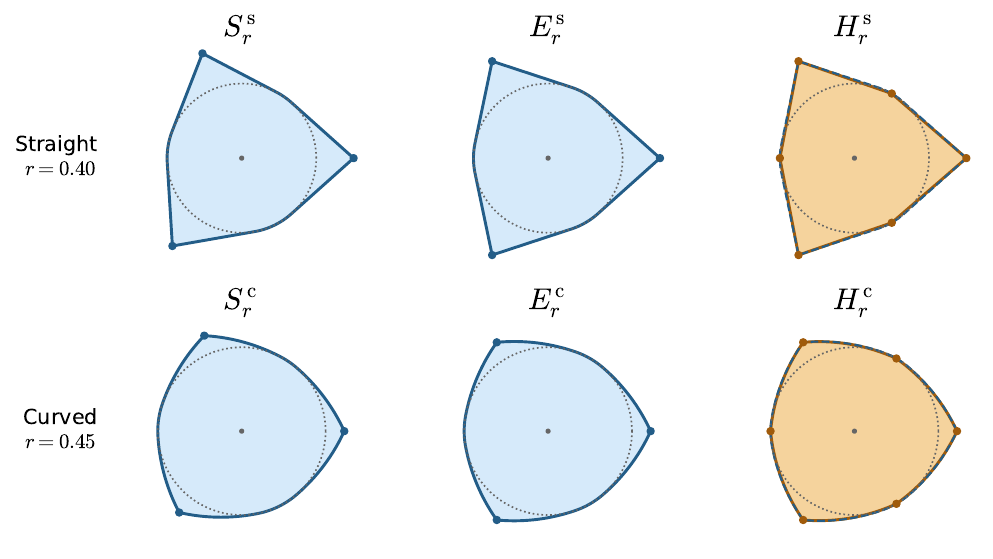}
\captionsetup{hypcap=false}
\captionof{figure}{The subsets in \eqref{eq:chain}--\eqref{eq:spectral-chain}.
Top: straight caps at $r=0.40$. Bottom: curved caps at $r=0.45$.
From left to right: an unequal cap body $S_r$, its equilateral
configuration $E_r$, and the inscribed hexagon $H_r$.
In the last column the dashed outline is $E_r$, so $H_r\subset E_r$;
the dotted circles have radius $r$. Angular balancing compares
functional values of $S_r$ and $E_r$, without asserting an inclusion.}
\label{fig:proof-overview}
\end{center}

\subsection{Map of the results}\label{sec:results-map}
Table~\ref{tab:pde-proof-map} identifies the steps proving the torsion
and first-eigenvalue assertions of Theorem~\ref{thm:main}.
The first reflection compares traces and normal derivatives on different
caps. The shape derivative requires a further comparison on the two
sides of the same cap; this is supplied by the sliding-source argument.
For the equilateral configurations, we study the inscribed hexagons
$H_r\subset E_r$ as $r$ varies. Thus this last step does not require
monotonicity of $E_r$ itself.

The geometric comparisons also give the area and Cheeger inequalities
through the inner-parallel-area theorem, Theorem~\ref{thm:inner}, and
\eqref{eq:cheeger-bounds}. Integrating these areas gives the distance-moment
bounds \eqref{eq:distance-moment}; Proposition~\ref{prop:distance-concavity}
proves their Minkowski concavity. Proposition~\ref{prop:counterexamples}
shows why concavity and inclusion monotonicity alone do not suffice for
the width inequalities.

For integral functionals, Theorem~\ref{kr:thm:clipping} and
Corollary~\ref{kr:cor:HL} give the radial comparisons.
Theorem~\ref{kr:thm:balancing} treats angular interaction energies;
Theorem~\ref{ck:thm:nonradial} and Corollary~\ref{kr:cor:small-s}
give the perturbative nonradial and Riesz results.

\begin{table}[H]
\small
\begin{center}
\begin{tabular}{|p{.3\textwidth}| p{.3\textwidth} | p{.3\textwidth} |}  
\hline 
{\bf Step} & {Minimal width: straight caps} & {Constant width: curved caps}\\
\hline 
Contained cap body\newline $S_r\subset K$
& Lemma~\ref{lem:straight-reduction}
& Lemma~\ref{lem:geometry}, \eqref{eq:SinK}\\ \hline 
Different caps: reflection comparing remote boundary data
& $T$: \eqref{eq:tracecomparison};\newline
  $\lambda_1$: \eqref{eq:eigen-remote}
& $T$: \eqref{eq:tracecomparison};\newline
  $\lambda_1$: \eqref{eq:eigen-remote}\\ \hline
Sliding boundary data into the attached half-cap
& $T$: Lemma~\ref{lem:slide}, \eqref{eq:straight-source-integral};\newline
  $\lambda_1$: Lemmas~\ref{lem:eigen-kernel} and~\ref{lem:eigen-cap-transfer}
& $T$: Lemma~\ref{lem:slide};\newline
  $\lambda_1$: Lemmas~\ref{lem:eigen-kernel} and~\ref{lem:eigen-cap-transfer}\\ \hline
Adjacent sides of one cap:\newline flux comparison
& $T$: \eqref{eq:cap-flux-comparison}--\eqref{eq:torsion-angular-derivative};\newline
  $\lambda_1$: \eqref{eq:eigen-angular-derivative}
& $T$: \eqref{eq:cap-flux-comparison}--\eqref{eq:torsion-angular-derivative};\newline
  $\lambda_1$: \eqref{eq:eigen-angular-derivative}\\ \hline
Equilateral optimizer: fixed $r$
& $T$: Theorem~\ref{prop:balance};\newline
  $\lambda_1$: Theorem~\ref{thm:eigen-balance}
& $T$: Theorem~\ref{prop:balance};\newline
  $\lambda_1$: Theorem~\ref{thm:eigen-balance}\\ \hline
Equilateral family:\newline inscribe $H_r\subset E_r$
& Construction \eqref{eq:straight-hexagon}
& Lemma~\ref{lem:hexgeom}\\ \hline 
Vary $r$ on $H_r$:\newline $T(H_r)\nearrow$, $\lambda_1(H_r)\searrow$
& $T$: Proposition~\ref{prop:hexmono};\newline
  $\lambda_1$: Proposition~\ref{prop:eigen-hexagon}
& $T$: Proposition~\ref{prop:hexmono};\newline
  $\lambda_1$: Proposition~\ref{prop:eigen-hexagon}\\ \hline
Complete the width inequalities and equality cases
& $T$: \S\ref{sec:torsion-final};\newline
  $\lambda_1$: \S\ref{sec:straight-disk-bounds}, \S\ref{sec:straight-final}
& $T$: \S\ref{sec:torsion-final};\newline
  $\lambda_1$: \S\ref{sec:curved-eigenvalue-final}\\
\hline 
\end{tabular}
\end{center}
\caption{Proof map for minimizing the torsional rigidity $T$ and maximizing $\lambda_1$.}
\label{tab:pde-proof-map}
\end{table}

\section{Cap subsets under width constraints}
\label{sec:cap-subsets}
Both width constraints admit subsets formed by a disk and three disjoint
caps. We describe the straight-cap and curved-cap constructions in turn.

\subsection{The straight-cap subset at minimal width}
\label{sec:straight-geometry}
We first describe the subset used for the minimal-width constraint.
Let $1/3\leq r<1/2$ and denote by
$\alpha=\arccos(r/(1-r))\in(0,\pi/3]$ the contact angle of a cap.
We consider the sets
\begin{equation}\label{eq:straight-caps}
 S_r^{\rm s}=\conv\bigl(\overline B_r,(1-r) e_{a_1},(1-r) e_{a_2},(1-r) e_{a_3}\bigr),
 \qquad \delta_i\geq2\alpha,\quad\sum_i\delta_i=2\pi,
\end{equation}
where $\delta_i$ are the consecutive angular gaps between the three
vertices. A cap with axis $a$ is bounded by two segments tangent to
$\overline B_r$ at $r e_{a\pm\alpha}$. Its radial function is
\begin{equation}\label{eq:straight-radial}
 p_a(\theta)=
 \begin{cases}
 r/\cos(\alpha-|\theta-a|),& |\theta-a|\leq\alpha,\\
 r,& |\theta-a|\geq\alpha,
 \end{cases}
\end{equation}
where angular distances are taken on the circle. The assumptions on
the gaps ensure that the three cap sectors have disjoint interiors.
The union of these caps and the disk has monotone boundary turning,
and is therefore convex. Thus this union is exactly the convex hull
in \eqref{eq:straight-caps}, without any additional segment joining
two caps.

\begin{lemma}\label{lem:straight-reduction}
A planar convex body of minimal width one has inradius
$1/3\leq r\leq1/2$. If $r<1/2$, it contains a straight three-cap
body \eqref{eq:straight-caps}, after translation. If $r=1/3$, the
original body is an equilateral triangle of altitude one.
\end{lemma}
\begin{proof}
We use the classical cap construction from
\cite[Chapter~6]{YaglomBoltyanskii}, also used in
\cite[Proposition~3.3]{LucardesiZucco} and \cite{widthpaper}.
The lower bound on $r$ is the planar case of Steinhagen's inequality
\cite{Steinhagen}. We recall its proof while constructing the subset.

We translate a largest inscribed disk so that its center is the origin.
The origin belongs to the convex hull of the contact normals.
Indeed, a strict separation would give a direction in which the disk
could be translated and enlarged. Since the disk is contained in the
body, we have $2r\leq1$. Suppose now that $r<1/2$. Two contact normals
cannot be antipodal, since the corresponding support lines would give
a width $2r<1$. We can therefore choose three contact normals
$n_1,n_2,n_3$ whose triangle contains the origin in its interior.
The corresponding support half-planes define a triangle $P\supset K$
circumscribed about $\overline B_r$.
If $a$ is its longest side and $p$ its perimeter, its least altitude
is $rp/a\leq3r$. Hence $1=w(K)\leq w(P)\leq3r$.

It remains to choose the three cap vertices. Let $V_i$ be the vertex opposite
the side with normal $n_i$. Since $h_K(n_i)=r$, the width condition
gives a point $x_i\in K$ satisfying $x_i\cdot n_i\leq-(1-r)$.
This point belongs to the corner region at $V_i$. To see this, observe
that each of the other corner caps is contained in
$\conv(\overline B_r,V_j)$, where $V_j\cdot n_i=r$.
Every point in such a cap satisfies $x\cdot n_i\geq-r>-(1-r)$.
We also have $|x_i|\geq1-r$, so convexity gives
$A_i=(1-r)x_i/|x_i|\in K$.
The cap $\conv(\overline B_r,A_i)\setminus\overline B_r$ lies in the
same corner region of $P$. Since the three corner regions have
disjoint interiors, the new cap sectors are disjoint as well.
Finally, convexity gives $\conv(\overline B_r,A_1,A_2,A_3)\subset K$.

If $r=1/3$, equality must hold in $p/a\leq3$, hence $P$ is
equilateral. For each of its three side normals, the opposite vertex
is the unique point of $P$ with support value $-2/3$.
The width condition therefore requires $K$ to contain all three
vertices, and we conclude that $K=P$.
\end{proof}

When $r=1/2$, the inclusion $\overline B_{1/2}\subset K$ will be
sufficient for the comparisons below. In the minimal-width class,
this inclusion does not imply that the original body is a disk.
We denote the equilateral straight-cap configuration by $E_r^{\rm s}$.
We will use these cap bodies as subsets; they need not have minimal
width one themselves.
\subsection{The curved-cap subset at constant width}\label{sec:geometry}
For constant width, we replace the convex hull used in the preceding
subsection by the unit-ball convex hull. We write $B_r$ for the open disk
of radius $r$ centered at the origin. If a set $X$ is contained in a
unit disk, its unit-ball convex hull is
\[
 \bh(X)=\bigcap\{\overline B_1(c):X\subset\overline B_1(c)\}.
\]
This set is convex. We refer to \cite{BallPolyhedra,ConstantWidthBook}
for the relation with spindle convexity. We first recall the two-point
case, which will allow us to identify circular arcs contained in a
unit-ball convex hull.

\begin{lemma}[Unit spindles]\label{lem:spindle}
For two distinct points at distance less than two, their unit-ball hull is the
closed spindle bounded by the two minor unit circular arcs through them.
\end{lemma}
\begin{proof}
After a rigid motion, we may take the two points to be
$(-a,0),(a,0)$, where $0<a<1$. Set $h=\sqrt{1-a^2}$ and consider
$L=\overline B_1((0,h))\cap\overline B_1((0,-h))$.
We show that every unit disk containing the endpoints also contains
this spindle. Its center $(s,t)$ satisfies $(|s|+a)^2+t^2\leq1$.
For any $(x,y)\in L$, we have $|x|\leq a$ and
$x^2+y^2-a^2+2h|y|\leq0$. We deduce that
\begin{align*}
 |(x,y)-(s,t)|^2
 &\leq1+x^2+y^2-a^2-2a|s|-2xs
       +2|y|\sqrt{1-(|s|+a)^2}\\
 &\leq1+x^2+y^2-a^2+2h|y|-2(a-|x|)|s|\leq1.
\end{align*}
Thus $L$ is contained in every admissible disk. Conversely, the two
disks defining $L$ contain both endpoints, so their intersection
contains the unit-ball convex hull. This proves the equality.
\end{proof}

\begin{lemma}\label{lem:geometry}
After translation, a constant-width-one body has concentric incircle and
circumcircle, with radii $r$ and $1-r$ satisfying
\[
 r_0\leq r\leq\tfrac12.
\]
If $r<1/2$, there are three circumcircle contacts $A_i=(1-r) e_{a_i}$ whose
cyclic angular gaps $\delta_i$ satisfy
\begin{equation}\label{eq:gaps}
 2\beta\leq\delta_i\leq\pi-\beta,\qquad
 \sum_i\delta_i=2\pi,\qquad
 \beta=\arccos\left(\frac1{2(1-r)^2}-1\right)\in(0,\pi/3].
\end{equation}
Moreover $K$ is an intersection of unit disks and contains
$\bh(\overline B_r\cup\{A_i\})$ for each $i$.
If $r=r_0$, then $K$ is a Reuleaux triangle of width one.
\end{lemma}
\begin{proof}
We first locate the centers of the inscribed and circumscribed disks.
For any center $c$, the support function relative to $c$ satisfies
$h_c(e)+h_c(-e)=1$, and therefore
\[
 \min_e h_c(e)=1-\max_e h_c(e),\qquad
 \max_e h_c(e)=\max_{x\in K}|x-c|.
\]
Thus a center minimizing the enclosing radius also maximizes the
inscribed radius. The center of the smallest enclosing circle belongs
to the convex hull of its contact points. Otherwise, moving the center
towards that convex hull would decrease the enclosing radius.
In the plane, at most three contact points are needed. If two of them
are antipodal, then $2(1-r)=1$, since $\operatorname{diam}K=1$.
In this case $\overline B_{1/2}\subset K\subset\overline B_{1/2}$,
and $K$ is the disk.

In the remaining case, we choose three contacts whose triangle contains
the center. Their angular gaps are at most $\pi$, and one of them is
at least $2\pi/3$. The diameter constraint then gives
$\sqrt3(1-r)\leq1$, which proves the asserted range of the radii.
This is the planar form of Jung's bound \cite{Jung}.
For each gap, $2(1-r)\sin(\delta_i/2)\leq1$ gives
$\delta_i\leq2\arcsin(1/(2(1-r)))=\pi-\beta$.
The lower bound for each gap follows by subtracting the two other
upper bounds from $2\pi$.

To prove the containment of the cap hulls, we use the identity
\[
 K=\bigcap_{y\in K}\overline B_1(y).
\]
The inclusion from left to right follows from the diameter constraint.
Conversely, if $x\notin K$, we can choose $e$ with
$x\cdot e>h_K(e)$ and then $y\in K$ with $y\cdot e=h_K(e)-1$.
Since $(x-y)\cdot e>1$, the point $x$ does not belong to
$\overline B_1(y)$. This proves the identity.
Every disk in this intersection contains $B_r$ and each $A_i$.
It therefore contains the corresponding unit-ball convex hull,
which proves the cap containment.

If $r=r_0$, the gap bounds force all three gaps to equal $2\pi/3$.
The contacts $A_1,A_2,A_3$ therefore form an equilateral triangle of
side one, and the diameter constraint gives
$K\subset\bigcap_i\overline B_1(A_i)=\mathcal R$.
Both bodies have constant width one. Their support functions are
ordered and have the same antipodal sums, so they coincide.
\end{proof}

\subsubsection{The explicit curved cap}
We now describe the cap at $A=(1-r) e_0$. Let
$L_A=\bh(\overline B_r\cup\{A\})$. We claim that its boundary outside
the disk consists of two unit circular arcs, with centers and upper
arc parametrization given by
\begin{equation}\label{eq:caparc}
 C_\pm=-(1-r) e_{\pm\beta},\qquad
 x(\theta)=C_++e_\theta,\quad \beta/2\leq\theta\leq\beta.
\end{equation}
The upper arc joins $A$ to $r e_\beta$, and the lower arc is its
reflection. The remaining boundary is $r e_\theta$ for
$\beta\leq\theta\leq2\pi-\beta$.
To prove the claim, let $\widetilde L$ be the region enclosed by this
curve. The outward normal angle is increasing, continuous at the
disk--cap joins and has a positive jump of $\beta$ at $A$.
Thus $\widetilde L$ is convex. By Lemma~\ref{lem:spindle}, each new
arc belongs to the spindle of its endpoints $A,r e_{\pm\beta}$.
The remaining circle is contained in $\overline B_r$. We obtain
$\widetilde L\subset\bh(B_r\cup\{A\})$ by convexity of the ball hull.

For the reverse inclusion, consider the disks
$D_\theta=\overline B_1(-(1-r) e_\theta)$ with
$\beta\leq\theta\leq2\pi-\beta$. Each contains $B_r$.
It also contains $A$, since this condition is equivalent to
$2(1-r)^2(1+\cos\theta)\leq1$, or $\cos\theta\leq\cos\beta$.
The disk is tangent at $r e_\theta$, and $D_{\pm\beta}$ contains
the corresponding exterior cap arc on its boundary.
Consequently, these disks give a supporting half-plane at every smooth
point of $\partial\widetilde L$. At $A$, their two extreme normals
generate its normal cone. Every separating normal at that corner is
a nonnegative combination of those two normals. Thus a point outside
$\widetilde L$ violates one of these supporting half-plane inequalities
and lies outside at least one $D_\theta$. We conclude that
\[
 \bh(B_r\cup\{A\})\subset\bigcap_\theta D_\theta
 \subset\widetilde L,
\]
which proves the description of $L_A$.

We will also need the radial function of the cap. The cap occupies
the polar angles $(-\beta,\beta)$, and for $0\leq\phi\leq\beta$ we have
\begin{equation}\label{eq:radialcap}
 p(\phi)=-(1-r)\cos(\beta-\phi)
       +\sqrt{1-(1-r)^2\sin^2(\beta-\phi)}.
\end{equation}
Since $\sqrt{1-(1-r)^2\sin^2 t}>(1-r)\cos t$, the derivative with respect to
$\beta-\phi$ is positive. Therefore $p$ decreases with $\phi$.
By reflection, the radial function decreases with the absolute angle
on both sides of the apex.

We can now combine the three caps by setting
\[
 S(r;a_1,a_2,a_3)=\bigcup_{i=1}^3 L_{A_i}.
\]
By \eqref{eq:gaps}, their polar sectors have disjoint interiors.
Along the counterclockwise boundary of the union, the outward normal
angle increases on every circular arc, is continuous at the disk--cap
joins and jumps forward by $\beta$ at each apex. The union is therefore
convex. Since each cap hull is contained in the original body, we have
\begin{equation}\label{eq:SinK}
 S(r;a_1,a_2,a_3)\subset K.
\end{equation}
We denote the equilateral configuration by
$E_r=S(r;0,2\pi/3,4\pi/3)$. As in the straight case, the comparison
will use the inclusion above; $S$ need not itself have constant width.

\begin{figure}[!htbp]
\centering\includegraphics[width=\linewidth]{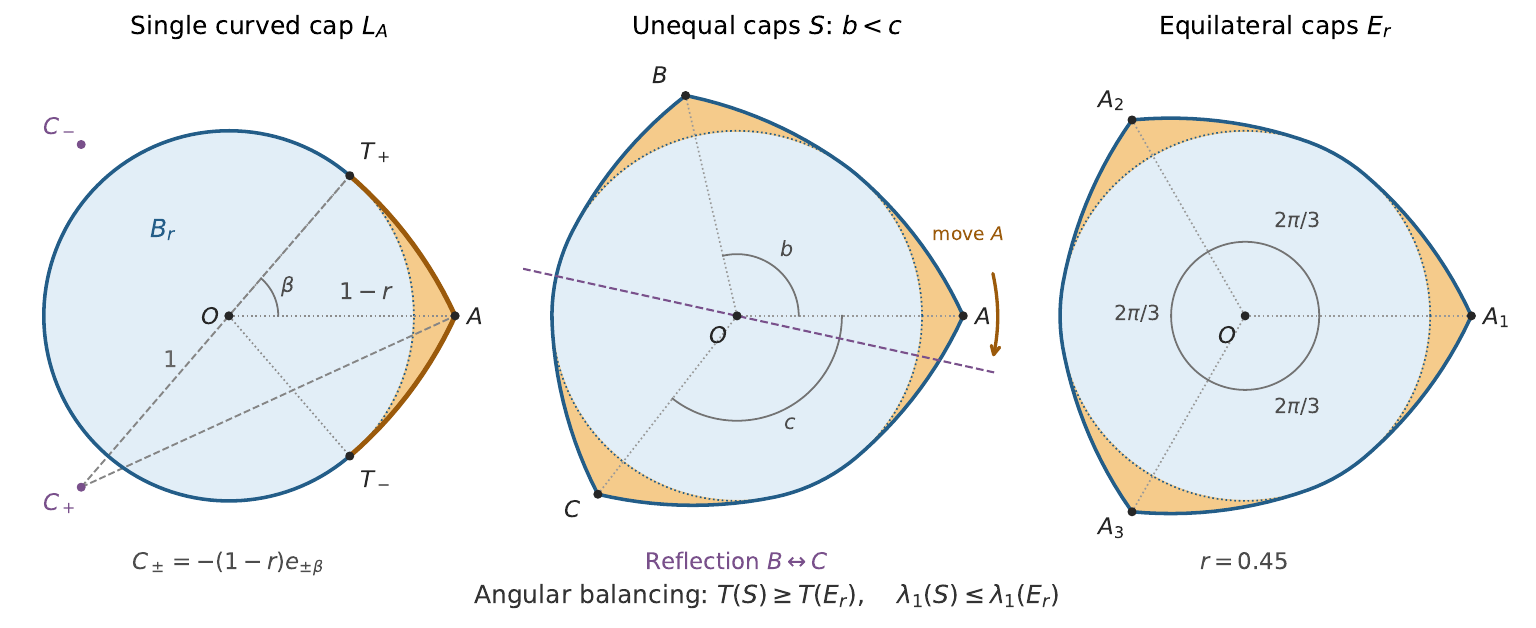}
\caption{The curved-cap construction at $r=0.45$. Left: a single cap,
its contact points $T_\pm$ with $B_r$, and the centers $C_\pm$ of its
unit circular arcs. Center: an unequal three-cap configuration with
$b<c$, the reflection interchanging the remote caps, and the clockwise
motion of $A$ that decreases torsion and increases the first eigenvalue.
Right: the equilateral comparison body $E_r$.}
\label{fig:caps}
\end{figure}
\subsection{Rotation of a single cap}
\label{sec:cap-rotation}
Fix the disk radius and rotate only the cap with axis zero
counterclockwise, keeping the other caps fixed. Let $\Gamma_+$ be its
open upper exterior branch and let $\sigma(x,y)=(x,-y)$, so that
$\Gamma_-=\sigma(\Gamma_+)$. This notation covers the straight segments
and the circular arcs. At smooth points of the moving boundary the
velocity is $Jx$, and its outward normal component is $V_n(x)=Jx\cdot n(x)$.
To determine its sign, write the single-cap radial boundary as
$x=p(\vartheta)e_\vartheta$. Both profiles above have $p'(\vartheta)<0$
on the open upper branch, and hence
\begin{equation}\label{eq:rotspeed}
 V_n(x)=-\frac{p(\vartheta)p'(\vartheta)}
                  {\sqrt{p(\vartheta)^2+p'(\vartheta)^2}}>0,
 \qquad V_n(\sigma x)=-V_n(x),\qquad x\in\Gamma_+.
\end{equation}
The second identity follows from $J\sigma=-\sigma J$ and
$n(\sigma x)=\sigma n(x)$. Reflection also preserves arclength.
The other caps and exposed incircle arcs have zero normal velocity.
The explicit speed formulas are recorded in
Appendices~\ref{app:variations} and~\ref{app:straight}; those appendices
also justify admissible one-sided rotations at cap contacts.
These paired speeds concern angular motion at fixed $r$. The hexagon
motions that vary $r$ will have different paired speeds.

\section{Inner parallel sets and Minkowski concavity}
\label{sec:geometric-inequalities}
We first derive the inner-parallel-area comparisons and their geometric
consequences, then examine the scope of inclusion monotonicity and
Minkowski concavity.

\subsection{Inner parallel areas, area and the Cheeger constant}
\label{sec:inner}
We use the inclusions proved above to compare inner parallel areas.
For a convex body $L$ and $t\geq0$, we denote
\[
 L_{-t}=\{x:x+\overline B_t\subset L\},\qquad A_L(t)=|L_{-t}|.
\]
We assign area zero to the empty set and to sets with empty interior.
Inner parallel sets preserve inclusion, and erosion commutes with
arbitrary intersections. These properties allow us to compute the inner
parallel sets of both cap families. We will show that, at a fixed depth,
their areas depend only on the original inradius, which reduces the area and Cheeger
comparisons to inequalities in one variable.

\subsubsection{Straight caps and homothetic erosion}
We first consider a straight three-cap body $S_r^{\rm s}$ from
Lemma~\ref{lem:straight-reduction}. We keep the notation
\[
 \alpha=\arccos\frac r{1-r},\qquad
 c(r)=\pi+3(\tan\alpha-\alpha),\qquad \frac13\leq r<\frac12.
\]
Each cap adds an area $r^2(\tan\alpha-\alpha)$ to the disk, and the
three added regions have disjoint interiors. To compute the inner
parallel sets, observe that the nonredundant supporting half-planes
have distance $r$ from the origin. Their normals correspond to the
remaining incircle arcs and to the tangent segments. If $G$ is the
complement of the three open angular cap sectors, we therefore have
$S_r^{\rm s}=\bigcap_{\theta\in G}\{x:x\cdot e_\theta\leq r\}$.
Moving each of these half-planes inward by $t$ gives
\begin{equation}\label{eq:straight-erosion}
 (S_r^{\rm s})_{-t}=(1-t/r)S_r^{\rm s},\qquad
 A_{S_r^{\rm s}}(t)=(r-t)^2c(r),\quad 0\leq t\leq r.
\end{equation}
At $t=r$, the inner parallel set reduces to the origin.
In particular, its area is independent of the positions of the caps.
Lucardesi and Zucco proved the following monotonicity at $t=0$ in
\cite[Lemma~3.2]{LucardesiZucco}. We use the homothety in
\eqref{eq:straight-erosion} to prove the corresponding statement at
positive depths.

\begin{lemma}\label{lem:straight-inner}
For each fixed $0\leq t<1/3$, the expression in
\eqref{eq:straight-erosion} is strictly increasing in
$r\in[1/3,1/2]$. It is positive at $t=1/3$ if $r>1/3$.
\end{lemma}
\begin{proof}
We differentiate with respect to $r$ on the open interval and obtain
\[
 \alpha'=-\frac1{(1-r)^2\sin\alpha},\qquad c'=3\tan^2\alpha\,\alpha'<0.
\]
Using $r=\cos\alpha/(1+\cos\alpha)$, we deduce that
\begin{align*}
 2c+rc'&=2\pi-6\alpha+3\tan\alpha(1-\cos\alpha)>0,\\
 \partial_r A_{S_r^{\rm s}}(t)
 &=(r-t)\{2c+rc'-tc'\}>0.
\end{align*}
The signs follow from $0<\alpha\leq\pi/3$ and $r>t$.
The formulas extend to the endpoints by continuity. In particular,
at $r=1/2$ the limiting body is $\overline B_{1/2}$ and $c=\pi$.
For $t=1/3$, positivity follows directly from
\eqref{eq:straight-erosion}.
\end{proof}

\subsubsection{Curved caps and erosion of their supporting disks}
For curved caps, we use an intersection of disks in place of the
supporting half-planes. Let $S_r^{\rm c}=S(r;a_1,a_2,a_3)$ be the body
constructed in Subsection~\ref{sec:geometry}, with cap angle
$\beta$. Denote by $G$ the complement of the three open sectors
$(a_i-\beta,a_i+\beta)$ on the angle circle. We have
\begin{equation}\label{eq:curved-disk-family}
 S_r^{\rm c}=\bigcap_{\theta\in G}\overline B_1(-(1-r) e_\theta),\qquad
 (S_r^{\rm c})_{-t}=\bigcap_{\theta\in G}\overline B_{1-t}(-(1-r) e_\theta)
 \quad(0\leq t\leq r).
\end{equation}
To prove the first identity, notice that the angle of each disk center
has distance at least $\beta$ from every apex axis. Thus each disk
contains the incircle and all three apex vertices. The disks corresponding to
sector endpoints contain the exterior cap arcs on their boundaries,
and the other disks are tangent along the remaining incircle arcs.
The supporting normal argument from Subsection~\ref{sec:geometry} then
gives equality. Since erosion commutes with intersection, it is enough
to reduce every disk radius by $t$ to obtain the second identity.

To describe the resulting boundary, set
\begin{equation}\label{eq:eroded-parameters}
 R=1-t,\qquad \rho=r-t=R-(1-r),\qquad
 \gamma=\arcsin\frac{(1-r)\sin\beta}{R},\qquad v=\beta-\gamma.
\end{equation}
For $0\leq t<r$, these parameters satisfy
$\beta/2\leq\gamma<\beta$ and $0<v\leq\beta/2$.
Consider a point $p e_\phi$, with $p\geq0$. Among the centers
$-(1-r) e_\theta$, $\theta\in G$, its distance is largest when the angle
in $G$ is nearest to $\phi$. On a remaining disk sector, the boundary
is therefore given by $p=\rho$. Inside a cap sector, the nearest
angle is one of the two endpoints of that sector. In particular, the
upper half of the cap with axis zero is parametrized by
\begin{equation}\label{eq:eroded-arc}
 x(\theta)=-(1-r) e_\beta+R e_\theta,\qquad \gamma\leq\theta\leq\beta.
\end{equation}
The first endpoint is the apex on the real axis, and the second is
$\rho e_\beta$. The three cap sectors remain disjoint for every
depth considered here. The same choice of nearest angles therefore
describes the active disks throughout the interval up to $r$.

We apply Green's area formula to these six arcs and obtain
\begin{equation}\label{eq:curved-inner-area}
 A_{S_r^{\rm c}}(t)=\pi\rho^2+
 3\{R^2v-(1-r)R\sin v-\rho^2\beta\},\qquad 0\leq t\leq r.
\end{equation}
Indeed, the integral of $x\times dx$ on the upper arc is
$R^2v-(1-r)R\sin v$, and the lower arc gives the same contribution.
Subtracting the area $\rho^2\beta$ of the replaced incircle sector
gives the contribution of one cap. As in the straight case, the
resulting area does not depend on the cap axes.

\begin{lemma}\label{lem:curved-inner}
For every fixed $0\leq t<r_0=1-1/\sqrt3$, the right side of
\eqref{eq:curved-inner-area} is strictly increasing in
$r\in[r_0,1/2]$. At $t=r_0$ it is positive when $r>r_0$.
\end{lemma}
\begin{proof}
It remains to determine the sign of the derivative with respect to $r$.
Since the area is independent of the cap axes, we may use the
equilateral configuration, which is defined throughout the inradius
interval. Differentiating $2(1-r)\cos(\beta/2)=1$ and using $(1-r)'=-1$ gives
\[
 \beta'=-\frac1{(1-r)^3\sin\beta},\qquad
 \frac1{(1-r)^2\sin\beta}=2\cot(\beta/2).
\]
For fixed $t$, the radius $R$ in \eqref{eq:eroded-arc} is constant,
while the derivative of the circle center is $e_\beta+2\cot(\beta/2)Je_\beta$.
The outward normal speed along the arc is therefore
$\cos(\beta-\theta)-2\cot(\beta/2)\sin(\beta-\theta)$.
The remaining incircle arcs have total angle $2\pi-6\beta$ and
outward normal speed one. Integrating the normal speeds along the
whole boundary, we obtain
\begin{equation}\label{eq:eroded-area-derivative}
 \partial_r A_{S_r^{\rm c}}(t)
 =\rho(2\pi-6\beta)+6R\{\sin v-2\cot(\beta/2)(1-\cos v)\}>0.
\end{equation}
One can also differentiate \eqref{eq:curved-inner-area} directly.
The endpoint terms cancel because the incident pieces share their
endpoints. At $r=r_0$, where the caps touch, the formula gives the
right derivative.

The first term in \eqref{eq:eroded-area-derivative} is nonnegative.
To prove that the second is positive, we use $0<v\leq\beta/2$ to write
\[
 \frac{\sin v}{1-\cos v}=\cot(v/2)
 \geq\cot(\beta/4)
 =\cot(\beta/2)+\csc(\beta/2)>2\cot(\beta/2).
\]
This proves the strict monotonicity. By continuity, the comparison
extends to the disk at $r=1/2$. Finally, if $t=r_0<r$, the inner
parallel set contains $B_{r-r_0}$ and therefore has positive area.
\end{proof}

For $t=0$, formula \eqref{eq:curved-inner-area} gives the three-cap
area bound used by Bezdek in \cite[Section~4]{BezdekDisk}:
\begin{equation}\label{eq:curved-area-zero}
 |S_r^{\rm c}|=\pi r^2+\frac32(1-2r^2)\beta
                    -\frac32\sqrt{4(1-r)^2-1}.
\end{equation}
The derivative in \eqref{eq:eroded-area-derivative} proves its
monotonicity and gives the same comparison at every erosion depth
in the lemma.

\subsubsection{Sharp comparisons and equality}
\begin{theorem}\label{thm:inner}
Let $K$ be a planar convex body of minimal width one, and let
$\mathcal T$ be the equilateral triangle of altitude one. Then
\[
 A_K(t)\geq A_{\mathcal T}(t)
 =\frac1{\sqrt3}(1-3t)_+^2\qquad(t\geq0).
\]
If instead $K$ has constant width one, then
\[
 A_K(t)\geq A_{\mathcal R}(t)\qquad(t\geq0).
\]
For either class, equality at any one $t\in[0,r(\mathcal Q)]$ forces
$K$ to be congruent to its corresponding optimizer $\mathcal Q$.
No equality characterization is asserted for $t>r(\mathcal Q)$.
\end{theorem}
\begin{proof}
Suppose first that the inradius $r$ is less than $1/2$. The geometric
construction gives $S_r\subset K$, so
$(S_r)_{-t}\subset K_{-t}$. We apply the preceding lemmas to compare
its area with that of the cap body at the smallest admissible inradius.
This body is $\mathcal T$ at $r=1/3$ in the straight case and
$\mathcal R$ at $r=r_0$ in the curved case. If $r=1/2$, we obtain
the same conclusion from the inscribed disk and the limiting formulas.
For depths beyond the inradius of the comparison triangle, its inner
parallel area is zero, so the asserted inequality also holds there.

We now consider equality in the stated interval. If $r>r(\mathcal Q)$,
the comparison is strict below $t=r(\mathcal Q)$, and the inner
parallel area remains positive at that endpoint. Consequently equality
requires $r=r(\mathcal Q)$. Lemma~\ref{lem:straight-reduction} gives
the conclusion in the straight case. In the curved case, $r=r_0$
forces the three gaps to equal $2\pi/3$, which gives
$\mathcal R\subset K$. Since the two bodies have constant width one,
their ordered support functions have the same antipodal sums.
They are therefore equal, and the bodies coincide.
\end{proof}

For $t=0$, we recover P\'al's and the Blaschke--Lebesgue area
inequalities. Ftouhi, Lucardesi and Saracco proved the minimal-width
inner-area comparison, including strictness, in
\cite[Corollary~4.5]{FtouhiLucardesiSaracco}. The constant-width
inner-area comparison and the stronger inner-perimeter comparison were
proved in \cite{BogoselCheeger}. The proof above derives both area
comparisons from the corresponding three-cap subsets.

We next apply these inequalities to the Cheeger constant. For a bounded
convex planar domain, it is defined by
\[
 h(K)=\inf_{E\subset K,\ |E|>0}\frac{P(E)}{|E|},
\]
where the competitors have finite perimeter. The Cheeger radius
$t_K=1/h(K)$ is the unique root in $(0,r(K))$ of
\begin{equation}\label{eq:cheeger-root}
 A_K(t_K)=\pi t_K^2;
\end{equation}
and the Cheeger set is $K_{-t_K}+\overline B_{t_K}$
\cite{KawohlLachand}. Applying the inner-area comparison at
$t_{\mathcal Q}$ gives $t_K\geq t_{\mathcal Q}$. We deduce that
\begin{equation}\label{eq:cheeger-bounds}
 \begin{aligned}
 h(K)&\leq h(\mathcal T)=3+\sqrt{\pi\sqrt3}
 &&\text{at minimal width one},\\
 h(K)&\leq h(\mathcal R)&&\text{at constant width one}.
 \end{aligned}
\end{equation}
The equality cases follow from Theorem~\ref{thm:inner}. For the Reuleaux
triangle, we can determine the sharp constant by substituting $r=r_0$
in \eqref{eq:curved-inner-area} and solving \eqref{eq:cheeger-root}.
This characterization does not require a numerical approximation.

\begin{remark}[Inner parallel perimeters]\label{rem:perimeters}
The same subset does not give the inner-perimeter comparison of
\cite{BogoselCheeger} through perimeter monotonicity. Indeed,
\[
 P(S_r^{\rm c})=2\pi r+3(1-2r)\beta
 =\pi+(1-2r)(3\beta-\pi)<\pi=P(\mathcal R)
 \quad(r_0<r<1/2).
\]
The perimeter of the subset is already smaller than the required bound
at $t=0$, and this remains true for small positive $t$ by continuity.
For minimal width, the triangle does not minimize perimeter in the
original class either, since
$P(B_{1/2})=\pi<2\sqrt3=P(\mathcal T)$.
Thus the inner-perimeter questions require a separate argument from
the inner-area comparisons proved above.
\end{remark}

\begin{figure}[tbp]
\centering
\begin{tabular}{@{}c@{\hspace{.02\linewidth}}c@{}}
\includegraphics[width=.48\linewidth]{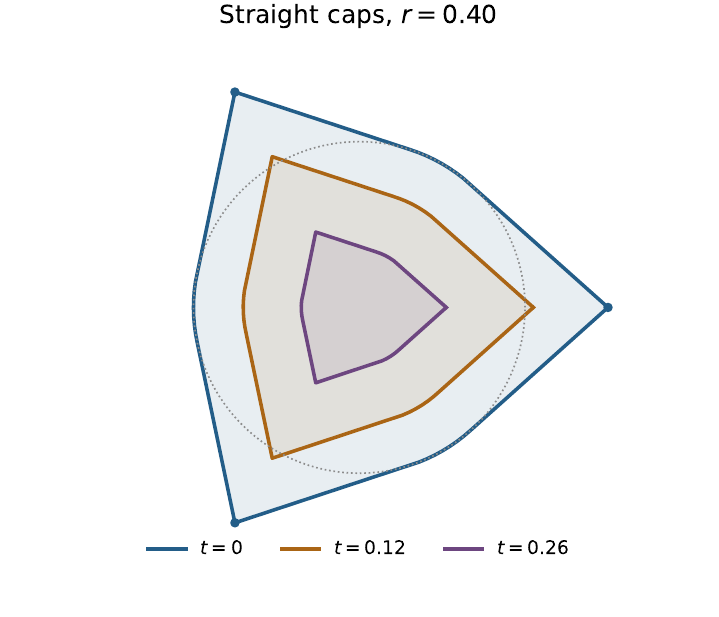} &
\includegraphics[width=.48\linewidth]{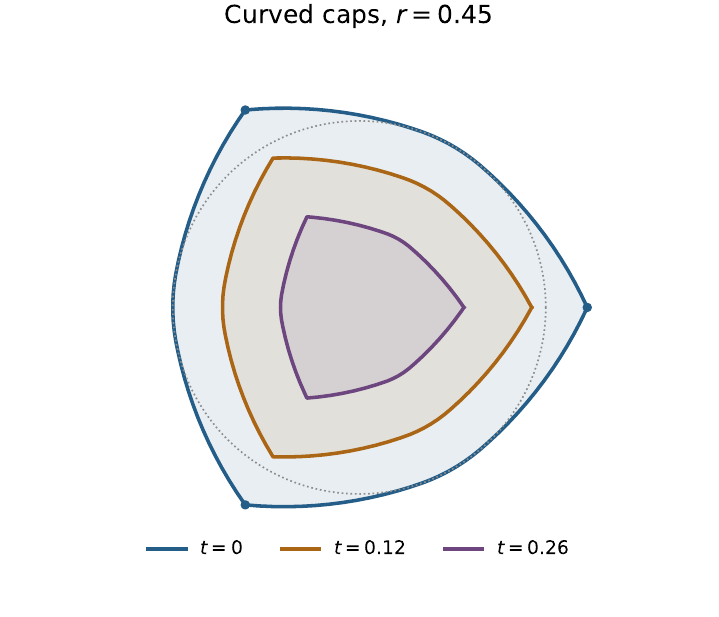}\\[1ex]
\includegraphics[width=.48\linewidth]{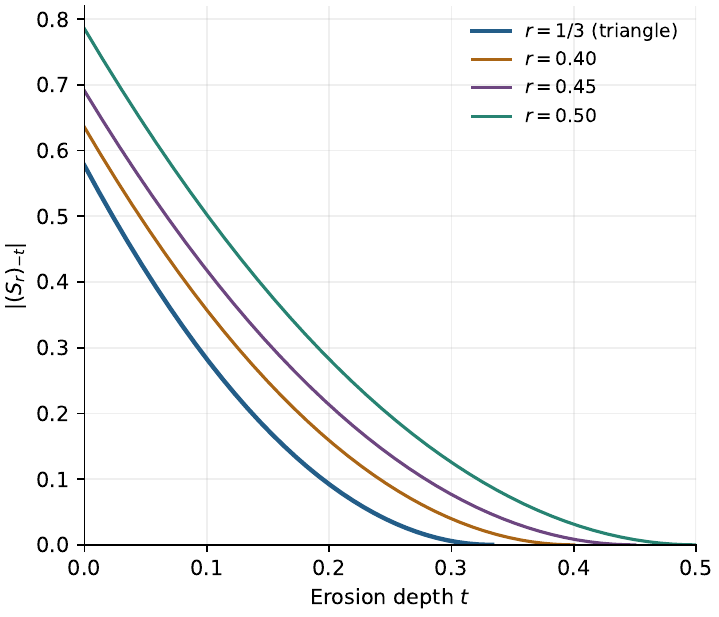} &
\includegraphics[width=.48\linewidth]{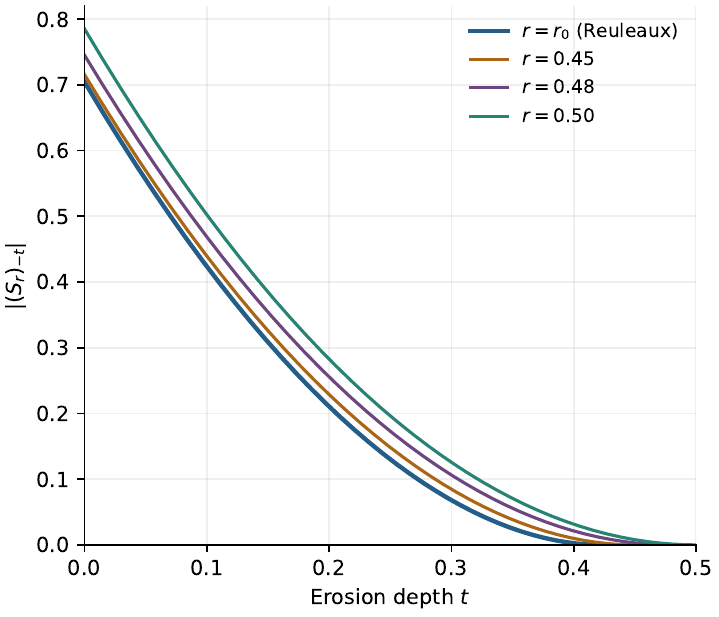}
\end{tabular}
\caption{Straight and curved three-cap bodies and their inner parallel
sets. Top: the straight family at $r=0.40$ and the curved family at
$r=0.45$, with erosion depths $t=0,0.12,0.26$.
Bottom: eroded areas for several values of the original inradius;
the endpoint curve is lowest in each class. Coordinates use the exact
radial formulas, sampled for rendering. The eroded bodies need not
have constant width.}
\label{fig:inner}
\end{figure}
\subsection{Minkowski concavity and further inequalities}
\label{sec:concavity}
We now examine whether the inclusion of the cap bodies can give
analogous comparisons for general concave functionals. We consider
Minkowski interpolation $K_\theta=(1-\theta)K+\theta L$, which differs
from taking the ordinary convex hull of $K\cup L$.
A functional $F$ is Minkowski concave if
\[
 F(K_\theta)\geq(1-\theta)F(K)+\theta F(L).
\]
The four problems considered here give the following functionals of
homogeneity degree one:
\begin{equation}\label{eq:bm-functionals}
 |K|^{1/2},\qquad h(K)^{-1},\qquad T(K)^{1/4},\qquad
 \lambda_1(K)^{-1/2}.
\end{equation}
They are increasing under inclusion, invariant under rigid motions
and Minkowski concave. For the PDE functionals, we refer to
\cite{bm_eigen,BorellTorsion,Colesanti}; for the Cheeger case see
\cite[Lemma~4.4]{FtouhiParallel} and the direct planar proof below.
We will show that these properties do not suffice to prove the
minimality of the equilateral or Reuleaux triangle. Concavity and
indecomposability have been used to study minimizers in
\cite{BucurFragalaLamboley}. For constant-width bodies, Kallay's
characterization \cite{Kallay} implies in particular that every
Reuleaux polygon is extreme. Knowing that a minimizer is extreme
therefore does not identify the triangle.

We first observe that Minkowski concavity does identify a maximizer
at constant width. If $F$ is continuous, invariant under rigid motions
and Minkowski concave, then
\begin{equation}\label{eq:concave-disk-maximum}
 F(K)\leq F(\overline B_{1/2})\qquad\text{if $K$ has constant width one}.
\end{equation}
Indeed, consider Minkowski averages of uniformly spaced rotations of $K$.
By concavity, the value of $F$ on each finite average is at least
$F(K)$. The support functions of these averages converge uniformly
to the constant $(2\pi)^{-1}\int h_K=1/2$. The averages therefore
converge to $\overline B_{1/2}$, and continuity gives
\eqref{eq:concave-disk-maximum}. This argument gives an upper bound
for $F(K)$, whereas the triangle comparisons require a lower bound.

\subsubsection{An inclusion criterion}
The following proposition gives a sufficient inclusion for such a
lower bound.
\begin{proposition}\label{prop:minkowski-certificate}
Let $F$ be increasing under inclusion, invariant under rigid motions,
and Minkowski concave. Fix a body $Q$. If
\begin{equation}\label{eq:certificate}
 \sum_{j=1}^N\theta_j g_jQ\subset K,
 \qquad\theta_j\geq0,\quad\sum_j\theta_j=1,
\end{equation}
where $g_j$ are rigid motions, then $F(K)\geq F(Q)$.
The same conclusion holds for Hausdorff limits of these combinations
when $F$ is continuous.
\end{proposition}
\begin{proof}
Using monotonicity and applying concavity to the finite sum, we obtain
$F(K)\geq F(\sum_j\theta_jg_jQ)\geq\sum_j\theta_jF(Q)=F(Q)$.
The statement for Hausdorff limits follows by continuity.
\end{proof}

Our geometric construction gives $S\subset K$. To apply the proposition
through this subset, we would also need a Minkowski average as in
\eqref{eq:certificate} contained in $S$. Perimeter shows that such an
inclusion is not always possible. For straight caps, take $K=\overline B_{1/2}$
and $Q=\mathcal T$. The inclusion would contradict the monotonicity
and Minkowski linearity of perimeter. For curved caps,
Remark~\ref{rem:perimeters} shows that no Minkowski average of
width-one Reuleaux triangles can be contained in $S_r^{\rm c}$ when
$r_0<r<1/2$.

There is another restriction at constant width. A Minkowski average
of bodies of constant width one still has constant width one.
Two such bodies cannot be properly nested: their ordered support
functions have the same antipodal sums and must therefore coincide.
Consequently, if $K$ and $Q$ both have constant width one, the inclusion
\eqref{eq:certificate} must be an equality. The next example shows
that an exact representation by averages of Reuleaux triangles is
not possible for every $K$.

\begin{proposition}[Counterexamples]
\label{prop:counterexamples}
There are continuous, positive, degree-one homogeneous,
rigid-motion invariant, inclusion-increasing, Minkowski-concave
functionals whose minimizer is not the proposed triangle in the
respective width class. In the constant-width class the counterexample
can also be chosen strictly increasing under proper inclusion.
\end{proposition}
\begin{proof}
For the minimal-width constraint, it is enough to consider perimeter,
since
$P(B_{1/2})=\pi<2\sqrt3=P(\mathcal T)$.
For the constant-width constraint, we define on planar convex bodies
\begin{equation}\label{eq:F5}
 F_5(K)=\min_{\theta\in[0,2\pi]}
 \int_0^{2\pi}h_K(\phi)\{1+\cos(5(\phi-\theta))\}\dd\phi.
\end{equation}
We first verify the properties of this functional. Since the density
is nonnegative and its first Fourier harmonic vanishes, each integral
is increasing under inclusion and invariant under translations.
The minimum over phases is invariant under rotations and reflections.
Moreover, each integral is Minkowski linear, so their infimum is
concave. Homogeneity of degree one and continuity for the Hausdorff
metric follow directly from the support-function formula.
To see positivity, we may translate an interior point to the origin.

We also obtain strict monotonicity. If $K\subsetneq L$, the support
difference $h_L-h_K$ is positive on an open angular interval.
Its integral against each density in \eqref{eq:F5} is therefore
strictly positive. These differences depend continuously on the phase,
so their minimum on the compact phase circle is positive.
It follows that $F_5$ is strictly increasing under proper inclusion.

To compare its values, we use $P(K)=\int_0^{2\pi}h_K$ and write
\[
 F_5(K)=P(K)-\left|\int_0^{2\pi}h_K(\phi)e^{5i\phi}\dd\phi\right|.
\]
If $\mathcal R$ is centered at its threefold symmetry center, its
fifth Fourier coefficient vanishes. Thus $F_5(\mathcal R)=\pi$.
Consider now the support function
\begin{equation}\label{eq:fivefold-support}
 h_{K_\varepsilon}(\phi)=\frac12+\varepsilon\cos5\phi,
 \qquad 0<\varepsilon<\frac1{48},
\end{equation}
whose antipodal support sum is one. Moreover,
$h+h''=1/2-24\varepsilon\cos5\phi>0$, so this defines a smooth,
strictly convex body of constant width one. A direct integration gives
\[
 F_5(K_\varepsilon)=\pi(1-\varepsilon)<\pi=F_5(\mathcal R).
\]
The Reuleaux triangle is therefore not a minimizer of $F_5$, which
completes the proof.
\end{proof}

The counterexamples can also be chosen strictly concave on
nonhomothetic bodies. Indeed, add $\delta|K|^{1/2}$ to either
functional, taking $\delta>0$ small enough to preserve the strict
comparison above. The strict Brunn--Minkowski inequality then gives
the required strict concavity \cite{Schneider}.
Consequently, for the functionals in \eqref{eq:bm-functionals},
the lower bounds require further properties of the functional or
of the comparison sets.

\subsubsection{Functionals defined by inner parallel areas}
We can obtain further sharp inequalities by integrating the comparison
in Theorem~\ref{thm:inner}. Let $\nu$ be a nonnegative Borel measure
on $[0,\infty)$ for which the following integrals are finite. Then
\begin{equation}\label{eq:weighted-erosion}
 \int_{[0,\infty)} A_K(t)\dd\nu(t)
 \geq \int_{[0,\infty)} A_{\mathcal Q}(t)\dd\nu(t),
\end{equation}
where $\mathcal Q$ is the triangle corresponding to the width class
under consideration. If $\nu([0,r(\mathcal Q)))>0$, equality holds
only for $\mathcal Q$. Indeed, every other body has strictly larger
inner parallel area at each depth in that interval. Thus the strict
comparison survives integration. This argument only uses the
inner-area inequality and requires no concavity assumption.

In particular, consider the moments of the distance to the boundary,
\begin{equation}\label{eq:distance-moment}
 M_p(K)=\int_K\operatorname{dist}(x,\partial K)^p\dd x
       =p\int_0^\infty t^{p-1}A_K(t)\dd t,\qquad p>0.
\end{equation}
We obtain $M_p(K)\geq M_p(\mathcal Q)$, with equality only for
$\mathcal Q$. For the equilateral triangle of minimal width one,
the integral can be computed explicitly:
\begin{equation}\label{eq:triangle-moment}
 M_p(\mathcal T)=\frac{2}{\sqrt3\,3^p(p+1)(p+2)}.
\end{equation}
For width $w$, the comparison scales by $w^{p+2}$.
These moments describe the distribution of distances to the boundary
and are different from torsional rigidity. The use of inner parallel
sets and distance-dependent test functions in estimates for torsion
is developed by Fragal\`a, Gazzola and Lamboley in
\cite{FragalaGazzolaLamboley}.

\begin{proposition}\label{prop:distance-concavity}
For every real $p>0$, $M_p^{1/(p+2)}$ is a degree-one,
inclusion-increasing, rigid-motion invariant Minkowski-concave
functional. Its sharp minimum in either width class is attained
uniquely by the corresponding triangle.
\end{proposition}
\begin{proof}
We prove the concavity by applying the Borell--Brascamp--Lieb
inequality to the distance functions. Write
$d_K(x)=\operatorname{dist}(x,\partial K)$ for $x\in K$, and extend
$f_K=d_K^p$ by zero outside $K$. Take $x\in\Int K$ and
$y\in\Int L$, and set $a=d_K(x)$ and $b=d_L(y)$.
The inclusions $x+\overline B_a\subset K$ and
$y+\overline B_b\subset L$ give
\[
 d_{K_\theta}((1-\theta)x+\theta y)\geq(1-\theta)a+\theta b.
\]
Raising this inequality to the power $p$, we obtain
\[
 f_{K_\theta}((1-\theta)x+\theta y)
 \geq\bigl((1-\theta)f_K(x)^{1/p}+\theta f_L(y)^{1/p}\bigr)^p
\]
whenever both input values are positive. In the
Borell--Brascamp--Lieb inequality, the power mean is taken to be zero
when either input is zero. Thus the displayed inequality is precisely
the hypothesis needed to apply that result. Taking parameter $1/p$
in dimension two, we deduce from \cite{BorellConvex,bm_eigen} that
\[
 M_p(K_\theta)^{1/(p+2)}\geq
 (1-\theta)M_p(K)^{1/(p+2)}+\theta M_p(L)^{1/(p+2)}.
\]
The remaining properties and the sharp comparison follow from
\eqref{eq:distance-moment} and Theorem~\ref{thm:inner}.
\end{proof}

When $m$ is a positive integer, we can also prove concavity directly
from the Brunn--Minkowski inequality in $\R^{2+m}$. For this purpose,
consider
\[
 \mathcal C_m(K)=\{(x,y):x\in K,\ |y|\leq\operatorname{dist}(x,\partial K)\}.
\]
On a convex body, the distance to the boundary is the infimum of the
affine distances to its supporting lines, and is therefore concave.
It follows that the set just defined is convex. Using the
characterization of distance by inscribed disks, we also have
\[
 (1-\theta)\mathcal C_m(K)+\theta\mathcal C_m(L)
 \subset\mathcal C_m((1-\theta)K+\theta L).
\]
Indeed, $x+B_a\subset K$ and $x'+B_b\subset L$ imply
$(1-\theta)x+\theta x'+B_{(1-\theta)a+\theta b}\subset K_\theta$,
and $|(1-\theta)y+\theta y'|\leq(1-\theta)|y|+\theta|y'|$.
Since $|\mathcal C_m(K)|=\omega_m M_m(K)$, applying the
Brunn--Minkowski inequality in dimension $m+2$ proves the desired
concavity.

Finally, we give a proof of inverse-Cheeger concavity in the plane
using inner parallel sets. Ftouhi states the result in arbitrary
dimension in \cite[Lemma~4.4]{FtouhiParallel}, by taking the $p\to1$
limit in the Brunn--Minkowski inequality of
\cite{ColesantiCuoghiSalani}. Let $a=1/h(K)$, $b=1/h(L)$ and
$s=(1-\theta)a+\theta b$. The definition of erosion gives
\[
 (1-\theta)K_{-a}+\theta L_{-b}\subset(K_\theta)_{-s}.
\]
The planar Brunn--Minkowski inequality, together with
\eqref{eq:cheeger-root}, implies
$|(K_\theta)_{-s}|^{1/2}\geq\sqrt\pi s$.
By uniqueness of the Cheeger root, we conclude that
$h(K_\theta)^{-1}\geq s$. As the preceding examples show, this
concavity property alone does not give the width-constrained
comparison; that conclusion uses the inner-area inequality.
\section{Torsional rigidity}
\label{sec:torsion}
We prove both torsion width inequalities by the same angular balancing
and paired-flux arguments. The two cap families have the common rotation
velocities of Subsection~\ref{sec:cap-rotation}; their inscribed hexagons
will then give the inradius comparisons.

The variational characterization shows that domain inclusion increases
$T$. This functional is continuous for Hausdorff convergence of convex
bodies with a uniform positive inradius. Indeed, uniform convergence of
the support functions gives inclusions between small homothetic
contractions and expansions of the limiting body. Continuity then
follows from domain monotonicity and scaling.

The convex domains considered here have boundaries consisting of
finitely many straight segments or circular arcs. Their torsion
functions belong to $H^2$ by convex-domain regularity
\cite[Chapter~3]{Grisvard}. On smooth boundary arcs they satisfy the
usual boundary regularity and Hopf lemma \cite{GT}. For a boundary
velocity $V$, the first-variation formula is
\begin{equation}\label{eq:hadamard}
 DT(\Omega)[V]=\int_{\partial\Omega}(\partial_nu_\Omega)^2V\cdot n\dd s.
\end{equation}
We use the shape differentiation results of
\cite[Sections~5.3 and~5.7]{HP}. For distributed and boundary shape
derivatives on nonsmooth domains, including curvilinear polygons, see
also \cite{Laurain}. Appendices~\ref{app:variations} and~\ref{app:straight}
give the pullback arguments needed when an admissible cap motion
separates a tangential contact. Thus \eqref{eq:hadamard} also holds for the one-sided
variations used below. At the endpoints of the interval for $r$, we
use continuity.

\subsection{Angular balancing for both cap families}\label{sec:torsion-angular}
The torsion comparisons follow the argument of
\cite[Section~3.1, Propositions~8--12]{widthpaper} for the minimal-width
problem. We do not reproduce all its details here. We recall the
boundary-kernel calculation and the interface comparison that allow the
same argument to cover both straight and curved caps: the transfer uses
the common circular attachment interface and does not require the exterior
cap boundary to be straight. Reflection first compares the two remote caps,
and the sliding-source argument then compares the two sides of the moving
cap. Figures~\ref{fig:remote-reflection} and~\ref{fig:same-cap-reflection}
illustrate these two reflections for both geometries.

\begin{lemma}[Moving a boundary source]\label{lem:slide}
Let $G$ be the upper half-disk of radius $r$ together with the upper half
of the cap at $A=(1-r) e_0$ from either family. Write $\gamma=\alpha$
for straight caps and $\gamma=\beta$ for curved caps. Boundary profiles are prescribed at
$r e_s$, $\gamma<s<\pi$, with zero Dirichlet data on the rest of $\partial G$.
If $f\geq0$ is a nonzero continuous boundary profile compactly supported in
$(\gamma,\pi)$ and $f_\tau(s)=f(s-\tau)$, where $\tau>0$ and both
supports remain in $(\gamma,\pi)$, then the harmonic extension of
$f_\tau$ is strictly smaller than that of $f$ inside the attached
half-cap. Its inward normal derivative is also strictly smaller on
the open exterior cap arc. The same conclusion holds for continuous
profiles vanishing at the endpoints when either support reaches
$\gamma$ or $\pi$.
\end{lemma}
\begin{figure}[!htbp]
\centering\includegraphics[width=\linewidth]{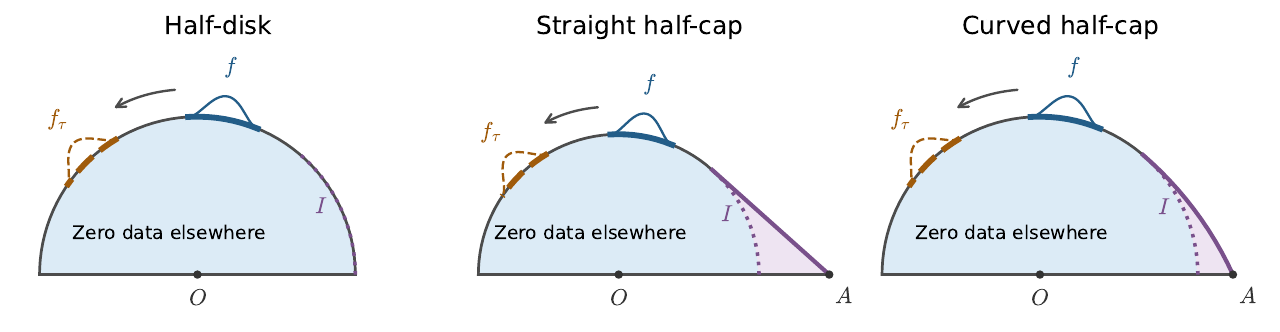}
\caption{Sliding boundary data in Lemma~\ref{lem:slide}.
The blue solid and ochre dashed profiles represent two separate Dirichlet
problems: the same nonnegative function is shifted toward larger angles,
with zero data on the remaining outer boundary. The small bumps indicate
boundary values, not changes to the domain. Left: the half-disk kernel
compares inward normal derivatives on the dotted arc $I$.
Middle and right: after attaching a straight or curved half-cap along $I$,
the later source gives smaller harmonic values in that half-cap and
smaller inward derivatives on its exterior boundary. The dotted arc is
then an interior interface. Boundaries use analytic segments and circles.}
\label{fig:sliding-source}
\end{figure}
\begin{proof}
We first take smooth boundary profiles, so that their harmonic extensions
belong to $H^1$. The construction is shown in
Figure~\ref{fig:sliding-source}. Scale the disk to radius one. The kernel mapping data at angle $s$ to
inward derivative at a disjoint zero-data point $t$ is
\[
 B_0(t,s)=\frac{\sin t\sin s}{\pi(\cos t-\cos s)^2},\qquad 0<t<s<\pi.
\]
This formula follows by differentiating the Poisson representation after
odd reflection across the diameter. Moreover,
\[
 \partial_s B_0(t,s)=\frac{\sin t}{\pi}
 \frac{\cos s\cos t+\cos^2s-2}{(\cos t-\cos s)^3}<0.
\]
Let $v$ be the solution with the earlier source and let $\psi$ be its
trace on the disk--cap interface $I$, whose angles are $(0,\gamma)$.
On the half-disk solve for $z$ with the same trace $\psi$ on $I$ and the
later source. The common interface data cancel when comparing normal
derivatives. Thus, with $n$ outward from the disk,
\[
 \partial_n z(e^{it})-\partial_n v(e^{it})
 =\int_{\operatorname{supp}f}f(s)
   \bigl[B_0(t,s)-B_0(t,s+\tau)\bigr]\,\mathrm ds>0,
 \qquad 0<t<\gamma.
\]
Paste $z$ on the disk to $v$ on the half-cap.
The resulting function is continuous and weakly superharmonic, since its $-\Delta$ is
the positive interface measure
$(\partial_n z-\partial_n v)\mathcal H^1|_I$.
It has the later source as its outer boundary data. The maximum principle
therefore bounds the later harmonic extension from above by this pasted
function. The positive interface measure gives strictness in the half-cap;
Hopf's lemma gives the derivative assertion on its open exterior arc.

The equality of the interface traces shows that the pasted function
belongs to $H^1$. No regularity of the isolated crescent at its tangential
tip is needed. For continuous profiles vanishing at their endpoints, we
use the Poisson extensions and approximate by smooth nonnegative data.
The off-diagonal kernel inequality gives the non-strict comparison in
the limit. Strictness follows by fixing a nonzero smooth subprofile
$0\leq g\leq f$ compactly supported away from the endpoints: the
comparison is strict for $g$ and nonnegative for $f-g$, and the
extensions depend linearly on the data. Hopf's lemma then gives the
strict derivative comparison. Thus the statement for continuous profiles
does not require their extensions to have finite energy.
\end{proof}

We will also use the corresponding point-source comparison. Let
$P_0(z,s)$ be the harmonic Poisson kernel on $G$, with angular source
measure on $\gamma<s<\pi$. The same interface argument gives
\[
 P_0(z,s)=\int_0^\gamma G_G(z,r e^{it})B_0(t,s)\,\mathrm dt,
 \qquad z\text{ in the attached half-cap},
\]
where $G_G$ is the Dirichlet Green function of $G$. Indeed, zero extension
of the half-disk solution with zero interface data has the negative
interface flux as its distributional $-\Delta$-source; Green representation
of the correction gives this formula. The normal-derivative and
arclength scale factors cancel. The positive Green weight and the strict
source monotonicity of $B_0$ show that $P_0(z,s)$ strictly decreases with
$s$. Its integral over the source arc is at most one, by the harmonic
maximum principle. Thus bounded truncated traces have absolutely
convergent Poisson representations. In addition, the explicit formula
for $B_0$ gives $P_0(z,s)=O(\pi-s)$ as $s\uparrow\pi$ at a fixed target.

\begin{theorem}[Angular balancing for torsion]\label{prop:balance}
At each fixed radius in the respective ranges
\[
 \begin{aligned}
 \frac13<r<\sqrt2-1&\quad\text{for the straight family }\eqref{eq:straight-caps},\\
 r_0<r<\frac12&\quad\text{for the curved family }\eqref{eq:gaps},
 \end{aligned}
\]
the equilateral cap body uniquely minimizes torsion, up to rotations.
\end{theorem}
\begin{proof}
Write $\gamma=\alpha$ for straight caps and $\gamma=\beta$ for curved caps.
Suppose that the gaps adjacent to an apex $A$ differ. After rotation,
write the axes as $0,b,-c$, with $2\gamma\leq b<c<\pi$ and
$b+c+2\gamma\leq2\pi$. These bounds follow from \eqref{eq:gaps} in the
curved case. In the stated straight-cap range, $\alpha>\pi/4$, so
\begin{equation}\label{eq:straight-small-gaps}
 2\alpha\leq\delta_i\leq2\pi-4\alpha<\pi.
\end{equation}
Let $S$ be the cap body and $u$ its torsion function.

Reflect across the line at angle $-(c-b)/2$, taking the side containing
axes zero and $b$ as positive. The disk and the remote caps are invariant.
The cap at zero is polarized toward the positive side, since its radial
profile decreases with angular distance from its axis. Thus the reflected
negative half of $S$ lies in its positive half. The difference of the two
torsion functions is harmonic and has nonnegative boundary values.
The maximum principle gives
\begin{equation}\label{eq:tracecomparison}
 u(r e_{b+t})\geq u(r e_{-c-t}),\qquad -\gamma<t<\gamma.
\end{equation}
The whole attachment arc at $b$ lies on that positive side: its angular
distance from the reflection line ranges between
$(b+c)/2-\gamma>0$ and $(b+c)/2+\gamma\leq\pi$.
At paired smooth exterior points of the remote caps, the same comparison
and Hopf's lemma order their inward fluxes. For moving the cap at zero,
we need the fluxes on its own two branches.

\begin{figure}[!htbp]
\centering\includegraphics[width=\linewidth]{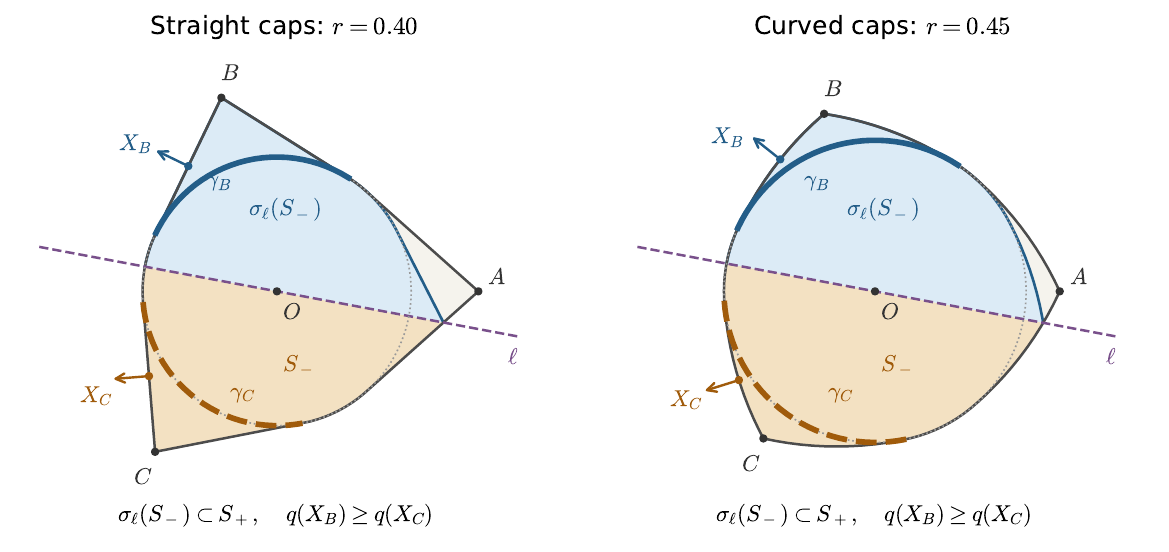}
\caption{The first reflection, comparing the remote caps at $B$ and $C$.
The dashed line $\ell$ bisects their axes. The ochre negative half $S_-$
reflects into the blue region $\sigma_\ell(S_-)\subset S_+$.
The highlighted attachment arcs $\gamma_B,\gamma_C$ are compared in
\eqref{eq:tracecomparison}; the same reflection orders the inward fluxes
$q=-\partial_nu$ at the paired exterior points $X_B,X_C$.
Arrows indicate outward unit normal directions, not flux magnitudes.
Left: straight caps at $r=0.40$. Right: curved caps at $r=0.45$.
Both configurations have axes $0,b,-c$, with $b=1.85<c=2.22$.
Boundaries are sampled from their analytic formulas.}
\label{fig:remote-reflection}
\end{figure}

Let $G$ be the half-disk with the upper half of the cap at zero attached,
as in Lemma~\ref{lem:slide}, and put $w=u-u\circ\sigma$ there, where
$\sigma(x,y)=(x,-y)$. This difference is harmonic, vanishes on the real
axis and exterior half-cap, and has semicircle data given by the source
at $b$ minus the reflected source at $-c$. If $c+\gamma\leq\pi$, both
remote caps lie in their respective half-planes. Inequality
\eqref{eq:tracecomparison} first orders their amplitudes; the
sliding-source comparison in Lemma~\ref{lem:slide} then gives $w>0$
in the attached half-cap.

\paragraph{The straight-cap crossing case.}
Only a straight remote cap can cross the negative real axis in the
ranges under consideration. If $c+\alpha>\pi$, define
\[
 a(s)=u(r e_{b-\alpha+s}),\qquad
 f(s)=u(r e_{-c+\alpha-s}),\qquad 0<s<2\alpha,
 \qquad m=\pi-c+\alpha.
\]
Here $0<m<2\alpha$, and \eqref{eq:tracecomparison} gives $a\geq f>0$.
On $0<s<m$, both source angles lie in $(\alpha,\pi)$, with
$b-\alpha+s<c-\alpha+s$; hence
$P_0(z,b-\alpha+s)a(s)>P_0(z,c-\alpha+s)f(s)$.
For $m<s<2\alpha$, the part of the cap at $-c$ lying above the axis
has angle $2\pi-c+\alpha-s\in(\alpha,\pi)$ and contributes positively.
The source decomposition is therefore
\begin{align}\label{eq:straight-source-integral}
 w(z)={}&\int_0^{2\alpha}P_0(z,b-\alpha+s)a(s)\dd s\notag\\
        &-\int_0^m P_0(z,c-\alpha+s)f(s)\dd s\notag\\
        &+\int_m^{2\alpha}P_0(z,2\pi-c+\alpha-s)f(s)\dd s>0.
\end{align}
The first two terms are strictly ordered on $(0,m)$, and the remaining
contributions are nonnegative. The preceding harmonic-kernel discussion
justifies splitting these bounded traces even at $s=\pi$.

In both geometries, Hopf's lemma on the open exterior half-cap now gives
\begin{equation}\label{eq:cap-flux-comparison}
 q_+(x)>q_-(x)>0,\qquad
 q_+(x)=-\partial_nu(x),\quad q_-(x)=-\partial_nu(\sigma x),
 \qquad x\in\Gamma_+.
\end{equation}
\begin{figure}[!htbp]
\centering\includegraphics[width=\linewidth]{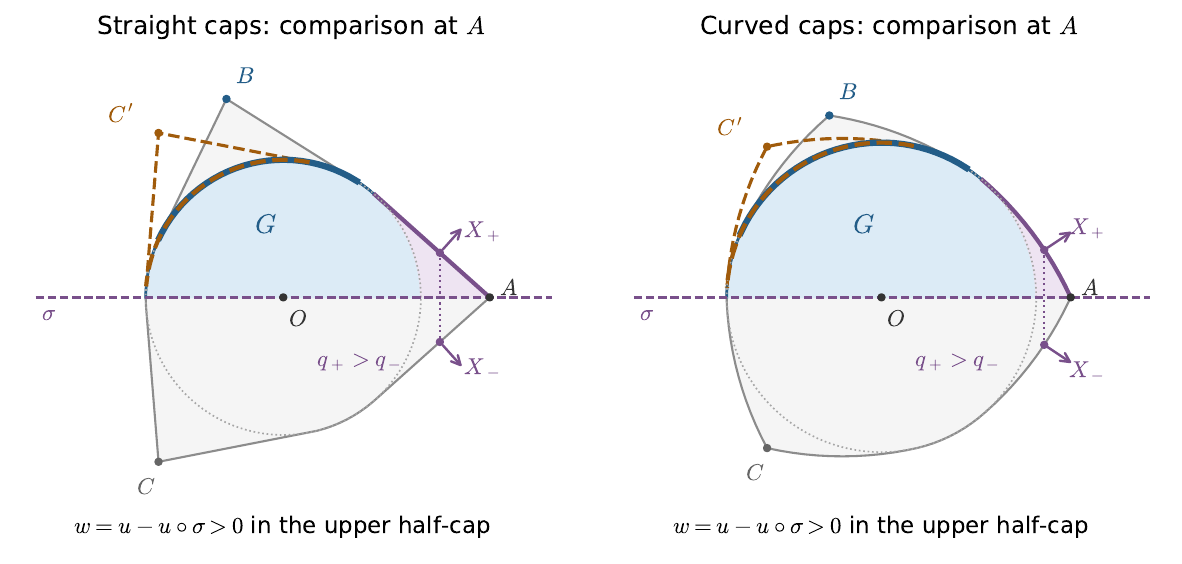}
\caption{The second reflection, comparing the two sides of the cap at $A$.
Here $\sigma(x,y)=(x,-y)$ and $C'=\sigma(C)$. The blue domain $G$ is
the upper half-disk with the upper half-cap at $A$ attached; the latter
is shaded violet. The dashed ochre cap and attachment arc are the images
of those at $C$. The first reflection compares the two boundary profiles;
Lemma~\ref{lem:slide} compares their angular positions. Together they give
$u-u\circ\sigma>0$ in the upper half-cap, and hence $q_+>q_-$ at its
paired exterior points $X_+,X_-$, as in \eqref{eq:cap-flux-comparison}.
Straight and curved configurations are the same as in
Figure~\ref{fig:remote-reflection}. Normal arrows indicate directions only.}
\label{fig:same-cap-reflection}
\end{figure}

By the paired velocities in \eqref{eq:rotspeed}, the Hadamard formula
has the same form for both cap families:
\begin{equation}\label{eq:torsion-angular-derivative}
 \frac{\partial T(S)}{\partial a_0}
 =\int_{\Gamma_+}V_n(x)(q_+^2-q_-^2)\dd s>0.
\end{equation}
A clockwise motion increases the smaller adjacent gap and decreases the
larger one, so is admissible and strictly decreases torsion.
At a straight-cap contact the speed is the same paired speed from
Subsection~\ref{sec:cap-rotation}; Appendix~\ref{app:straight} justifies
its one-sided derivative. Appendix~\ref{app:variations} gives the
corresponding justification for curved caps. There is no extra endpoint
term in either case.

The admissible gap sets are compact, and torsion depends continuously
on the gaps. A minimum therefore exists. At every unequal gap triple,
the preceding admissible variation decreases torsion, so precisely the
equilateral configurations minimize it.
\end{proof}

\subsection{The inscribed hexagons}\label{sec:hexagon-geometry}
The equilateral cap bodies contain hexagons with alternating vertex radii.
The straight family uses ordinary edges; the curved family uses unit
circular arcs. Both constructions are available on their full parameter
intervals.

\paragraph{Straight hexagons.}
For the straight family we use the hexagon
\begin{equation}\label{eq:straight-hexagon}
 H_r^{\rm s}=\conv\{(1-r) e_{2j\pi/3},r e_{(2j+1)\pi/3}:j=0,1,2\},
 \quad \frac13\leq r\leq\frac12.
\end{equation}
Its vertices alternate between cap vertices and incircle points, so
$H_r^{\rm s}\subset E_r^{\rm s}$. It is a convex hexagon for $r>1/3$;
at $r=1/3$ the radius-$r$ vertices are the three side midpoints of
$\mathcal T$. Its symmetries are $D_3$.

\paragraph{Curved hexagons.}
For $r\in[r_0,1/2]$, take the six successive vertices
\[
 v_{2j}=(1-r) e_{2j\pi/3},\qquad
 v_{2j+1}=r e_{(2j+1)\pi/3},\qquad j=0,1,2.
\]
We replace each side of this hexagon by its outward minor unit circular
arc and denote the enclosed body by $H_r$. This family will allow us to
compare the equilateral cap configurations for different values of $r$.

\begin{lemma}\label{lem:hexgeom}
The body $H_r$ is a convex intersection of six unit disks (three distinct
disks at $r=r_0$), has dihedral symmetry $D_3$, and satisfies
\[
 H_r\subset E_r,\qquad H_{r_0}=\mathcal R.
\]
\end{lemma}
\begin{proof}
Consider $A=(1-r,0)$ and $B=(r/2,\sqrt3r/2)$. Define
\begin{equation}\label{eq:center}
 \begin{gathered}
 L=|B-A|=\sqrt{1-3r+3r^2},\qquad h=\sqrt{1-L^2/4},\\
 n_m=\frac{(\sqrt3r/2,\,1-r-r/2)}L,\qquad
 C=\frac{A+B}{2}-h n_m.
 \end{gathered}
\end{equation}
The required circle has center $C$ and its minor arc has angle
$\alpha=2\arcsin(L/2)$. The other centers are obtained by the reflections
and rotations of $D_3$.

We verify that these disks support the proposed boundary. Write $C=(X,Y)$.
The disk contains $A,B$ on its boundary. We have $Y<0$ and
$3X+\sqrt3Y=3/2-\sqrt3h/L<0$.
For the first inequality, squaring its two positive sides reduces it to
$(4-L^2)(2-3r)^2-3r^2L^2=12(1-r)^2(1-3r^2)>0$.
The second follows from $L^2<1$. For the other four vertices, their
squared distances to $C$ minus one are, respectively,
\[
 (1-r)(3X-\sqrt3Y),\quad r(3X+\sqrt3Y),\quad
 (1-r)(3X+\sqrt3Y),\quad 2\sqrt3rY.
\]
Thus the only remaining vertex comparison is
\[
 |(1-r) e_{2\pi/3}-C|^2-1=(1-r)(3X-\sqrt3Y)\leq0.
\]
It is equivalent to
$h(3r-1)/L\geq\sqrt3(1-r)/2$. Both sides are positive and its squared
form is
\begin{equation}\label{eq:vertexalgebra}
 (4-L^2)(3r-1)^2-3(1-r)^2L^2
 =-12r^2(3r^2-6r+2)\geq0.
\end{equation}
Two points on a unit circle that lie in another distinct unit disk have
their minor circle arc in that disk: the allowed arc on the circle has
length less than $\pi$. Hence all six proposed arcs satisfy all six disk
constraints. Writing $\theta_m=\arg n_m$, the endpoint normal angles
are $\theta_m\pm\alpha/2$. Their order within $[0,\pi/3]$ follows from
\[
 \theta_m-\alpha/2\geq0\ \Longleftrightarrow\ 3r^2\leq1,
 \qquad
 \theta_m+\alpha/2\leq\pi/3
 \ \Longleftrightarrow\ 3r^2-6r+2\leq0.
\]
Both inequalities hold on the prescribed interval. It follows from the order of these normals that the six arcs form the
whole boundary of the convex disk intersection.

Each outer arc between $v_i,v_{i+1}$ belongs to the spindle
$\bh(\{v_i,v_{i+1}\})$. One endpoint is an apex $A_j$ and the other
belongs to $B_r$. By monotonicity of ball hulls this spindle lies in
$L_{A_j}\subset E_r$. Since $E_r$ is convex, it contains the region
enclosed by the six arcs. At $r_0$, $L^2=2-\sqrt3$ and $\alpha=\pi/6$;
successive pairs of circle centers coincide with the three Reuleaux
vertices. The six arcs therefore join to the three boundary arcs of
$\mathcal R$.
\end{proof}

\begin{figure}[tbp]
\centering
\begin{tabular}{@{}c@{\hspace{.02\linewidth}}c@{}}
\includegraphics[width=.48\linewidth]{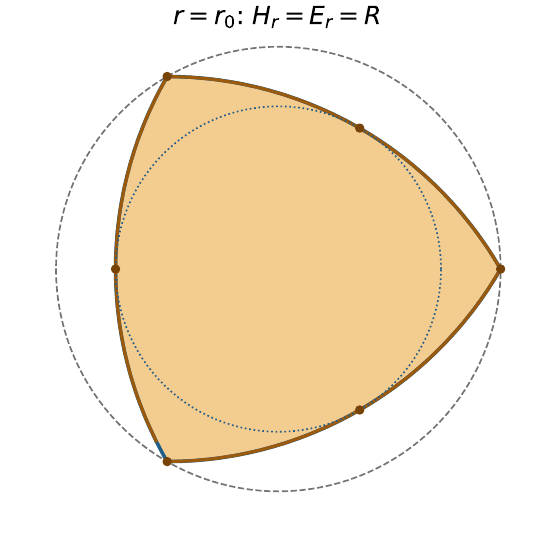} &
\includegraphics[width=.48\linewidth]{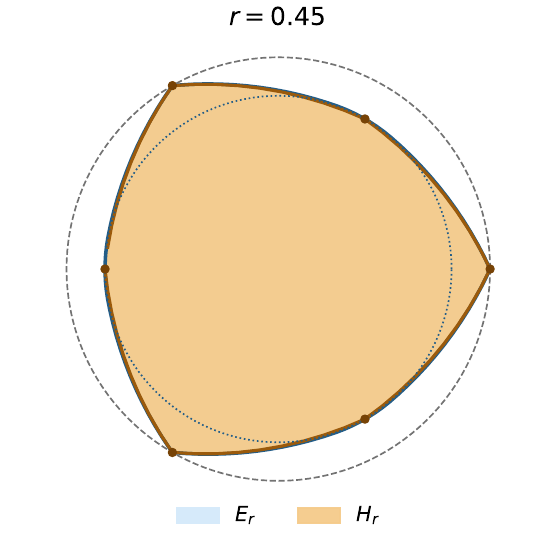}\\[1ex]
\includegraphics[width=.48\linewidth]{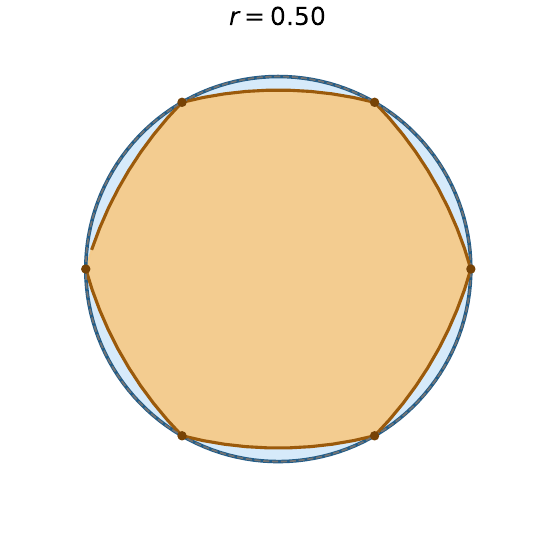} &
\includegraphics[width=.48\linewidth]{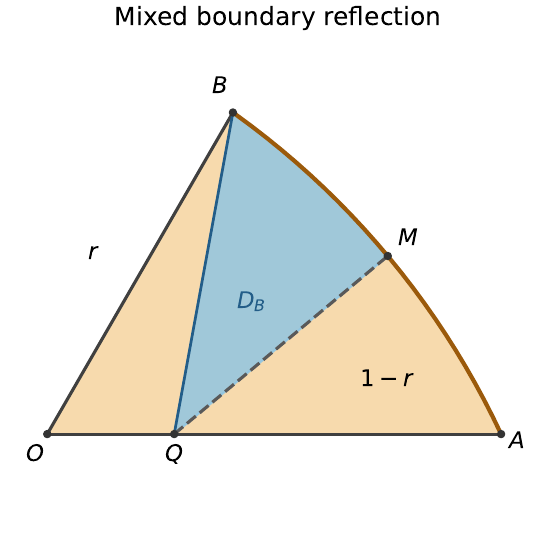}
\end{tabular}
\caption{The equilateral curved-cap body $E_r$ (blue) and its disk hexagon
$H_r$ (ochre). Dotted and dashed circles have radii $r$ and $1-r$.
At $r_0$ the two bodies coincide with the Reuleaux triangle. At $r=1/2$
the cap body is the disk, while the hexagon remains a proper subset.
The lower-right panel shows the reflected sector subdomain $D_B$ at $r=0.45$.
All boundaries come from analytic circular arcs; plotted outlines are sampled.}
\label{fig:comparison}
\end{figure}

\subsection{Monotonicity of torsion on the hexagons}\label{sec:torsion-hexagon}
To compare the boundary fluxes for both hexagon families, we first
determine a cone containing the gradient of the torsion function. We use only symmetry
and convexity of the superlevel sets. Concavity, reflection and gradient cones for more
general sectors with mixed boundary conditions are studied in
\cite{BucurFragalaMixed}. In our case the sector has angle $\pi/3$ and
is obtained by restricting a solution on either symmetric convex hexagon.
\begin{lemma}[A gradient cone]\label{lem:gradient-cone}
Let $H$ be either $H_r^{\rm s}$ or $H_r$ from
Subsection~\ref{sec:hexagon-geometry}. Let $v\in C^1(\Int H)$ have convex
superlevel sets and be invariant under the symmetries of $H$. In the sector $0<\arg x<\pi/3$,
\begin{equation}\label{eq:gradientcone}
 \nabla v=-c_0 e_0-c_1 e_{\pi/3},\qquad c_0,c_1\geq0.
\end{equation}
\end{lemma}
\begin{proof}
Each superlevel set contains the segment joining a point to its reflection
in either bounding ray. Differentiating at that point toward its reflected
point gives $v_y\leq0$ for the first ray and
$-\sqrt3v_x+v_y\geq0$ for the second. These inequalities are equivalent
to \eqref{eq:gradientcone}, with
$c_1=-2v_y/\sqrt3$ and $c_0=-v_x+v_y/\sqrt3$.
\end{proof}

\begin{lemma}[Paired boundary fluxes]\label{lem:flux}
Let $H$ be either hexagon, with $r$ in the interior of its parameter
interval, and let $u$ be its torsion function. On the boundary piece from $A=(1-r) e_0$ to $B=r e_{\pi/3}$, points paired by reflection
across its perpendicular bisector satisfy
\[
 q_B\geq q_A,
\]
where the subscript identifies the half closer to the indicated endpoint.
\end{lemma}
\begin{proof}
The sector between the rays $OA,OB$ has Dirichlet data on its outer
arc or segment and homogeneous Neumann data on both rays. The torsion
function on the full domain has convex superlevel sets \cite{MakarLimanov}; this
also follows from concavity of $\sqrt u$ on convex domains
\cite{Kennington}. These concavity properties, as well as the
log-concavity used below for the first eigenfunction, extend to the
convex bodies considered here by approximation with smooth convex
bodies and interior convergence of the solutions. Uniqueness gives their reflection
symmetries. Lemma~\ref{lem:gradient-cone} therefore applies to $u$.

Let $M$ be the midpoint of the outer arc or segment. The perpendicular
bisector of $AB$, passing through $M$, meets $OA$ at
\[
 Q=(\xi,0),\qquad \xi=\frac{(1-r)^2-r^2}{2(1-r)-r}=\frac{1-2r}{2-3r}.
\]
Since $0<\xi<1-r$, the region bounded by $QA,QB$ and the outer piece
$AB$ is contained in the sector and invariant under reflection $\sigma$
in $QM$. Its half $D_B$ adjoining $B$ has boundary piece $BM$, median
$MQ$, and segment $QB$. The reflected segment $QB$ is $QA$, where the
normal derivative of $u$ vanishes.

Let $w=u-u\circ\sigma$ on $D_B$. Then $\Delta w=0$, $w=0$ on $BM\cup MQ$,
and $\partial_n w=\partial_nu\geq0$ on $QB$, with $n$ pointing out of
$D_B$. Indeed, an unnormalized such normal is
$(-\sqrt3r/2,r/2-\xi)$; its scalar product with
$-c_0 e_0-c_1 e_{\pi/3}$ equals
$\sqrt3(c_0r+c_1\xi)/2\geq0$.
The mixed maximum principle, obtained by testing with $w_-$, gives
$w\geq0$. On the open arc or segment $BM$, the boundary point comparison gives
$\partial_n u\leq\partial_n(u\circ\sigma)$, and hence $q_B\geq q_A$.
This is the flux comparison needed in the following proposition.
\end{proof}

\begin{proposition}\label{prop:hexmono}\label{prop:straight-hexagon}
The functions $r\mapsto T(H_r^{\rm s})$ on $[1/3,1/2]$ and
$r\mapsto T(H_r)$ on $[r_0,1/2]$ are strictly increasing.
\end{proposition}
\begin{proof}
For the curved family, fix $r_0<r<1/2$ and let $t_m=Jn_m$. We differentiate the circle
center in \eqref{eq:center}, keeping its radius equal to one. We obtain
\begin{equation}\label{eq:velocity}
 C'=a_r n_m+b_r t_m,\qquad
 a_r=(1-2r)\left(\frac{\sqrt3}{4L}-\frac3{8h}\right)>0,\qquad
 b_r=\frac1{4L}+\frac{\sqrt3h}{2L^2}>0.
\end{equation}
In this computation both the chord midpoint and its normal depend on
$r$. In particular $\theta_m'=-\sqrt3/(2L^2)$ and
$h'=3(1-2r)/(8h)$. Positivity of $a_r$ follows from
$2h>\sqrt3L$, equivalent to $L^2<1$.

At paired arc parameters $\theta_m\pm s$, $0<s<\alpha/2$, the normal
speeds are $a_r\cos s\pm b_r\sin s$. The plus sign corresponds to the
shorter-radius endpoint $B$. By symmetry and the Hadamard formula,
\begin{align}\label{eq:positive}
 \frac{\dd}{\dd r}T(H_r)
 =6\int_0^{\alpha/2}\bigl[
 a_r\cos s(q_B^2+q_A^2)+b_r\sin s(q_B^2-q_A^2)\bigr]\dd s>0.
\end{align}
The second term is nonnegative by Lemma~\ref{lem:flux}; the first is
strictly positive on the open arc. The assertion at the two endpoints
follows by continuity.
At $r=1/2$ the derivative is zero by the symmetry interchanging the two
alternating radii; this does not affect strict increase on the interval.

For the straight family, work in the sector with $A=(1-r)e_0$ and
$B=r e_{\pi/3}$. Put $L=|B-A|=\sqrt{1-3r+3r^2}$.
Lemma~\ref{lem:flux} gives the same paired flux order $q_B\geq q_A>0$.
The outward normal to $AB$ is
$n=(\sqrt3r/2,1-r-r/2)/L$. As $r$ increases, its normal speed $b$ is
affine along the segment, with
\[
 b(A)=-\frac{\sqrt3r}{2L},\qquad
 b(B)=\frac{\sqrt3(1-r)}{2L}.
\]
For paired points $U$ in the $B$-half and $V$ in the $A$-half,
\[
 b(U)>0,\qquad b(U)+b(V)=\frac{\sqrt3(1-2r)}{2L}>0.
\]
Thus
\[
 q_B^2b(U)+q_A^2b(V)
 =(q_B^2-q_A^2)b(U)+q_A^2\{b(U)+b(V)\}>0.
\]
Applying the torsion Hadamard formula on the full domain and summing
over its six congruent edges gives $T'(H_r^{\rm s})>0$ for
$1/3<r<1/2$. These variations are justified by piecewise affine maps of the six
sector triangles. By continuity at the endpoint polygons, the result
holds on the closed interval.
\end{proof}

\subsubsection{Proof of the torsion width inequalities}
\label{sec:torsion-final}\label{sec:curved-torsion-final}
\label{sec:straight-torsion-final}
For minimal width, the inscribed disk gives $T(K)\geq\pi r^4/8$.
Introduce the torsion cutoff
\begin{equation}\label{eq:torsion-cutoff}
 r_T=\left(\frac{2\sqrt3}{45\pi}\right)^{1/4}<\frac25<\sqrt2-1.
\end{equation}
The bound follows from $\pi>25/8$ and $\sqrt3<7/4$.
If $r>r_T$, the disk estimate is strictly larger than $T(\mathcal T)$.
At $r=r_T$, strictness follows because $K$ properly contains its
inscribed disk: its minimal width is one, while $2r_T<1$.
The strong maximum principle gives strict torsion monotonicity under
this inclusion. For $1/3<r<r_T$, the straight-cap reduction, angular
balancing in Theorem~\ref{prop:balance}, and
Proposition~\ref{prop:hexmono} give
\[
 T(K)\geq T(S_r^{\rm s})\geq T(E_r^{\rm s})
 \geq T(H_r^{\rm s})>T(H_{1/3}^{\rm s})=T(\mathcal T).
\]
At $r=1/3$, Lemma~\ref{lem:straight-reduction} identifies $K$ with
$\mathcal T$.

For constant width and $r_0<r<1/2$, Lemmas~\ref{lem:geometry} and
\ref{lem:hexgeom}, Theorem~\ref{prop:balance}, and
Proposition~\ref{prop:hexmono} give \eqref{eq:chain}, with a strict
final inequality. At $r=1/2$, $K$ is the disk and contains $H_{1/2}$,
so the same strict conclusion holds. At $r=r_0$, the gap bounds force
all three gaps to equal $2\pi/3$, and $\mathcal R=S\subset K$.
Two constant-width-one bodies cannot be properly nested: their ordered
support functions have equal antipodal sums. Thus $K=\mathcal R$.
Scaling proves both torsion assertions and their equality cases in
Theorem~\ref{thm:main}.

\section{Dirichlet Laplace eigenvalues}
\label{sec:eigenvalues}
We prove the first-eigenvalue comparison in both width classes. The
cap reflections, exceptional straight-cap source decomposition, and
paired velocities have been established in Sections~\ref{sec:cap-subsets}
and~\ref{sec:torsion}. We now supply the spectral gaps and Helmholtz
kernel comparisons required to use that geometry for eigenfunctions.

The Rayleigh characterization shows that domain inclusion decreases
$\lambda_1$. As for torsion, homothetic inclusions and the scaling
$\lambda_1(a\Omega)=a^{-2}\lambda_1(\Omega)$ give Hausdorff continuity
on convex bodies with a uniform positive inradius. On the convex,
piecewise smooth domains considered here, the positive first
eigenfunction belongs to $H^2$ \cite[Chapter~3]{Grisvard} and satisfies
the usual boundary regularity and Hopf lemma on smooth arcs \cite{GT}.
For a boundary velocity $V$, its first-variation formula is
\begin{equation}\label{eq:eigen-hadamard}
 D\lambda_1(\Omega)[V]=-\int_{\partial\Omega}(\partial_n\varphi)^2
                  V\cdot n\dd s,\qquad\int_\Omega\varphi^2=1.
\end{equation}
The shape differentiation results and the one-sided pullback argument
recalled in Section~\ref{sec:torsion} apply equally to
\eqref{eq:eigen-hadamard}; at endpoint configurations we again use
continuity.

\subsection{Helmholtz kernels and angular balancing for both cap families}
\label{sec:eigen-angular}

The two reflections used for torsion remain the geometric starting point.
For torsion their difference is harmonic, whereas reflected copies of
a first eigenfunction give a difference satisfying
$(-\Delta-\lambda)w=0$, with $\lambda=\lambda_1(\Omega)$.
Thus each comparison requires a spectral gap on its auxiliary domain.
There is a second issue: the boundary data are a positive profile minus
another positive profile at a different angular position. They need not
have a pointwise sign. To compare them we replace the explicit harmonic
kernel of Lemma~\ref{lem:slide} by a Helmholtz kernel, and prove that
moving a positive source to a larger angle decreases its influence.
We first work on the half-disk and then transfer this property through
the disk--cap interface. This transfer applies to both cap geometries.
The spectral kernel order will follow from a positive separated Green
expansion, replacing the explicit Poisson-kernel calculation used for torsion.

We identify $\R^2$ with $\mathbb C$ and write $j_{0,1}$ for the first
positive zero of $J_0$. If $G$ is connected and
$0\leq\lambda<\lambda_1(G)$, then
$-\Delta-\lambda$ satisfies the Dirichlet maximum principle on $G$.
Indeed, testing a weak supersolution with its negative part and using
the Rayleigh principle makes that negative part vanish. A nonzero
nonnegative solution is strictly positive in the interior, and the Hopf
boundary lemma applies on smooth boundary pieces. In the angular
comparisons below, the gap follows from inclusion in $\Omega$, with a
nonempty open region omitted. Indeed, equality of the two Dirichlet
eigenvalues would make a zero-extended ground state a ground state on
the larger connected domain, contradicting positivity on that omitted
region. The mixed boundary comparison in Subsection~\ref{sec:eigen-hexagon}
will require a separate spectral-gap argument.

For the kernel comparison, fix
$U=\{\rho e^{it}:0<\rho<1,\ 0<t<\pi\}$ and $0\leq\kappa<j_{0,1}$,
so that $\kappa^2<\lambda_1(U)$. The Poisson kernel $P_\kappa(z,s)$
is characterized by the following representation: for
$f\in C_c^\infty(0,\pi)$, the unique solution of the Dirichlet problem
\[
 \begin{cases}
  (-\Delta-\kappa^2)h_f=0&\text{in }U,\\
  h_f(e^{is})=f(s)&0<s<\pi,\\
  h_f=0&\text{on }(-1,1)
 \end{cases}
 \qquad\text{is}\qquad
 h_f(z)=\int_0^\pi P_\kappa(z,s)f(s)\,\mathrm ds.
\]
Here $\mathrm ds$ is angular measure, equal to arclength on the unit semicircle.
Equivalently,
\[
 P_\kappa(z,s)=-\left.\partial_{n_y}G_{U,\kappa}(z,y)\right|_{y=e^{is}},
\]
where $n_y$ is outward and $G_{U,\kappa}$ is the Dirichlet Green function
normalized by $(-\Delta_y-\kappa^2)G_{U,\kappa}(z,y)=\delta_z$.
For Green representation and the variational Dirichlet theory, see
\cite[Sections~2.4 and~8.2]{GT}.

\begin{lemma}[Helmholtz boundary kernel]\label{lem:eigen-kernel}
With $U$, $\kappa$ and $P_\kappa$ as above, define, for $t\ne s$,
\[
 B_\kappa(t,s)=\lim_{\rho\uparrow1}
       \frac{P_\kappa(\rho e^{it},s)}{1-\rho}.
\]
Here $s$ is the source angle and $t$ is a target angle where the boundary
data vanish. Thus $B_\kappa$ measures the inward normal response at $t$
to a boundary source at $s$; it is the double inward normal derivative
of the Dirichlet Green kernel at these distinct boundary points.
For $0<t<s<s'<\pi$,
\begin{equation}\label{eq:eigen-kernel}
 0<B_\kappa(t,s')\leq
 \frac{\pi-s'}{\pi-s}B_\kappa(t,s)<B_\kappa(t,s).
\end{equation}
In particular, moving the same nonnegative boundary profile to larger
source angles strictly decreases its inward normal response on every
earlier zero-data arc.
\end{lemma}

\begin{proof}
We use the change of variables $x=-\log\rho$ to separate the radial
and angular variables. After multiplication by $\rho^2$, the operator becomes
\[
 L_\kappa+A,
 \qquad L_\kappa=-\frac{\mathrm d^2}{\mathrm dx^2}-\kappa^2e^{-2x},
 \qquad A=-\frac{\mathrm d^2}{\mathrm dt^2},
\]
on $(0,\infty)\times(0,\pi)$ with Dirichlet boundary conditions.
The Green function is unchanged by this transformation: the conformal
factor cancels the Jacobian of the point source. At $\rho=1$ the inward normal derivative
is $\partial_x$.

The bound $\kappa<j_{0,1}$ is the unit-disk spectral threshold,
$\kappa^2<\lambda_1(B_1)=j_{0,1}^2$; it also gives a positive function
$h(x)=J_0(\kappa e^{-x})$ on $[0,\infty)$ satisfying $L_\kappa h=0$.
This is the step that makes the separated expansion useful. For a finite
strip length $M>0$, integration by parts gives the ground-state identity
\begin{equation}\label{eq:eigen-groundstate}
 \int_0^M\bigl(|f'|^2-\kappa^2e^{-2x}f^2\bigr)\,\mathrm dx
 =\int_0^Mh^2\bigl|(f/h)'\bigr|^2\,\mathrm dx
 \quad (f\in H^1_0(0,M)).
\end{equation}
On each finite interval $h$ is bounded above and away from zero.
Poincar\'e's inequality applied to $f/h$ therefore makes the Dirichlet
eigenvalues $\mu_{j,M}$ of $L_\kappa$ on $(0,M)$ strictly positive.
Write $\psi_{j,M}$ for the corresponding real orthonormal eigenfunctions.
The Green kernel of $A+\mu$ on $(0,\pi)$ is, for $t<s$,
\[
 g_\mu(t,s)=
 \frac{\sinh(\sqrt\mu\,t)\sinh(\sqrt\mu(\pi-s))}
 {\sqrt\mu\sinh(\sqrt\mu\pi)}.
\]
At $\mu=0$ this expression is understood by continuity as
$t(\pi-s)/\pi$. On the rectangle $(0,M)\times(0,\pi)$, expand the
Green kernel in the $x$-eigenfunctions and take the two inward normal
derivatives at the boundary $x=0$. The result is
\begin{equation}\label{eq:eigen-positive-series}
 B_{\kappa,M}(t,s)=
 \sum_{j=1}^\infty |\psi'_{j,M}(0)|^2g_{\mu_{j,M}}(t,s),
 \qquad t\ne s.
\end{equation}
For fixed $M$, differentiation is justified away from the diagonal. The min--max
principle gives
$(j\pi/M)^2-\kappa^2\leq\mu_{j,M}\leq(j\pi/M)^2$; the one-dimensional
ODE and trace estimates bound $|\psi'_{j,M}(0)|$ by a polynomial in
$1+\mu_{j,M}$. On compact sets where $|t-s|>0$, exponential decay in
$\sqrt\mu|t-s|$ dominates those polynomial bounds. The series and the
needed local derivatives therefore converge uniformly there.

We now compare the angular factors. For $0<a'<a$ and $k\geq0$,
$\sinh(ka')/\sinh(ka)\leq a'/a$, with the continuous interpretation
at $k=0$. Consequently each term of \eqref{eq:eigen-positive-series} satisfies
the first inequality in \eqref{eq:eigen-kernel}. Its coefficient is
$|\psi'_{j,M}(0)|^2\geq0$, because the same boundary derivative occurs
in both variables. Consequently summation preserves the comparison.

It remains to justify the limit of these differentiated kernels.
From \eqref{eq:eigen-groundstate}, first for compactly supported functions,
and the angular Poincar\'e inequality, the half-strip form satisfies
\[
 \langle(L_\kappa+A)v,v\rangle\geq\|v\|_{L^2}^2.
\]
The bounded potential then also gives coercivity in $H^1_0$.
Zero extension identifies the rectangle form domains with increasing
subspaces whose union is dense in the half-strip space. Coercivity gives
strong $H^1$ convergence of their variational solutions for fixed
$L^2$ right-hand sides. Domain comparison makes the positive rectangle
Green kernels $K_M$ increase. Testing against smooth compactly supported
sources identifies their limit with the half-strip Green kernel $K$;
local elliptic regularity gives pointwise convergence off the diagonal.

Fix $0<t<s<\pi$ and choose disjoint boundary half-balls $V_t,V_s$
about $(0,t),(0,s)$, away from the horizontal sides. On their product,
$W_M:=K-K_M$ is nonnegative, bounded by the smooth off-diagonal kernel $K$,
and tends to zero in $L^2(V_t\times V_s)$. It solves the homogeneous
equation in each variable and has zero trace on each flat boundary.
For fixed $Y\in V_s$, a local $L^2$-to-$C^1$ boundary estimate in the
first variable gives
\[
 |\partial_{n_1}W_M((0,t),Y)|
 \leq C_t\|W_M(\cdot,Y)\|_{L^2(V_t)}.
\]
Set $a_M(Y):=\partial_{n_1}W_M((0,t),Y)$. Since the neighborhoods
are separated and away from the diagonal and corners, joint local
elliptic regularity permits this differentiation; $a_M$ satisfies
the same homogeneous equation in the second variable on $V_s$ and
has zero trace on its flat Dirichlet boundary. The same boundary
estimate in that variable, followed by integration of the first bound,
therefore gives
\[
 |\partial_{n_2}a_M((0,s))|\leq C_s\|a_M\|_{L^2(V_s)},
 \qquad
 \|a_M\|_{L^2(V_s)}\leq C_t\|W_M\|_{L^2(V_t\times V_s)}.
\]
Consequently,
\[
 \left|\partial_{n_1}\partial_{n_2}W_M((0,t),(0,s))\right|
 \leq C_{t,s}\|W_M\|_{L^2(V_t\times V_s)}\longrightarrow0,
\]
where both normals point inward. The two neighborhoods are fixed and
separated, lie on flat Dirichlet boundary pieces away from the corners,
and have coefficients independent of $M$. Odd reflection of the
zero-trace functions, with even extension of the potential, reduces
each boundary estimate to a local interior elliptic gradient estimate
\cite{GT}. Thus $C_{t,s}$ is independent of large $M$, uniformly on
compact sets of separated boundary pairs.
Here $B_\kappa(t,s)=\partial_{n_1}\partial_{n_2}K((0,t),(0,s))$,
with the same formula for $B_{\kappa,M}$ using $K_M$: since
$x=-\log\rho\sim1-\rho$ as $\rho\uparrow1$, this agrees with the
defining limit $\lim_{\rho\uparrow1}P_\kappa(\rho e^{it},s)/(1-\rho)$,
so the mixed-derivative convergence gives
$B_{\kappa,M}(t,s)\to B_\kappa(t,s)$ off the diagonal.
Under the logarithmic map the rectangles correspond to half-annuli
exhausting $U$. Their Dirichlet form domains exhaust $H^1_0(U)$,
so this limit imposes no additional boundary condition at the origin.

Passing to the limit gives the weak inequality in
\eqref{eq:eigen-kernel}. The Poisson kernel is positive because
$\kappa^2<\lambda_1(B_1)<\lambda_1(U)$. It vanishes on a boundary neighborhood of $t$
disjoint from its source, so Hopf's lemma gives $B_\kappa(t,s)>0$.
The factor $(\pi-s')/(\pi-s)<1$ proves the strict inequality.
\end{proof}

\begin{remark}[Sliding a positive boundary profile]
	Lemma~\ref{lem:eigen-kernel} gives a comparison for nonnegative integrable boundary
	profiles. Let \(0<\alpha<\pi\), let
	\(f\in L^1(\alpha,\pi)\), \(f\ge 0\), \(f\not\equiv0\), be supported in
	\((\alpha,\pi)\), extend it by zero outside this interval, and for
	\(\tau>0\) set
	\[
	f_\tau(s)=f(s-\tau),
	\]
	assuming that both \(\operatorname{supp} f\) and
	\(\operatorname{supp} f_\tau\) are contained in \((\alpha,\pi)\).
	In particular, \(\operatorname{supp}f\subset(\alpha,\pi-\tau)\).
	If \(q_f(t)\) denotes the inward normal response on the zero-data arc
	\(0<t<\alpha\), then the separated source and target give the finite integral
	\[
	q_f(t)
	=
	\int_\alpha^\pi B_\kappa(t,s)f(s)\,\mathrm ds .
	\]
	After the change of variables \(s=u+\tau\), the support conditions give
	\[
	q_{f_\tau}(t)
	=
	\int_\alpha^{\pi-\tau} B_\kappa(t,u+\tau)f(u)\,\mathrm du .
	\]
	For \(u\in\operatorname{supp}f\), we have \(t<\alpha<u<u+\tau<\pi\),
	so Lemma~\ref{lem:eigen-kernel} yields
	\[
	B_\kappa(t,u+\tau)<B_\kappa(t,u).
	\]
	Since \(f\) is nonnegative and nonzero in \(L^1\), integration gives
	\[
	q_{f_\tau}(t)<q_f(t),
	\qquad 0<t<\alpha .
	\]
	Thus translating a fixed nonnegative boundary source toward the endpoint
	\(\pi\) strictly decreases its inward normal response on every earlier
	zero-data arc. Under its geometric and spectral hypotheses,
	Lemma~\ref{lem:eigen-cap-transfer} shows that this ordering is preserved after the
	response is transmitted through the circular interface into an attached
	region.
\end{remark}

\begin{lemma}[Transfer through the circular interface]
\label{lem:eigen-cap-transfer}
Let $0<\alpha<\pi$ and
$\Gamma=\{e^{it}:0<t<\alpha\}$. Let $D\supset U$ be a bounded,
connected Lipschitz domain. Define the attached open region by
\[
 F:=D\setminus\overline U,
 \qquad D=U\cup\Gamma\cup F.
\]
Assume that $F$ is connected and
$\partial U\setminus\overline\Gamma\subset\partial D$, so the
attachment is through $\Gamma$.
Assume
$0\leq \kappa<j_{0,1}$ and $\kappa^2<\lambda_1(D)$.
For the Poisson kernel $P_{D,\kappa}$ of $-\Delta-\kappa^2$,
with source on the remaining semicircle and zero Dirichlet data
on the rest of the exterior boundary,
\[
 P_{D,\kappa}(z,s)>P_{D,\kappa}(z,s')
 \quad(z\in F,\ \alpha<s<s'<\pi).
\]
The corresponding inward normal derivative kernels satisfy the same
strict inequality on smooth outer boundary pieces of $F$ away from
$\overline\Gamma$. Both the straight and circular half-caps used below
satisfy these geometric hypotheses. Lipschitz regularity is required
of $D$, not of the isolated region $F$, which may have a tangential cusp.
\end{lemma}

\begin{proof}
Let $G_{D,\kappa}$ be the Dirichlet Green function on $D$. We express
the Poisson kernel on $D$ in terms of the half-disk kernel by proving
\begin{equation}\label{eq:eigen-transfer}
 P_{D,\kappa}(z,s)=\int_0^\alpha
 G_{D,\kappa}(z,e^{it})B_\kappa(t,s)\,\mathrm dt,
 \qquad z\in F,\quad \alpha<s<\pi.
\end{equation}
The boundary source first induces a positive inward flux density on
the circular interface; the positive Green kernel then propagates
this interface source into the attached region.
To this end, take $f\in C_c^\infty((\alpha,\pi))$ and let $v$ solve the equation on
$U$ with data $f$ on the remaining circle and zero elsewhere, and let
$\varphi$ solve the corresponding problem on $D$. Set $V=v$ on $U$
and $V=0$ on $F$; the pieces match on $\Gamma$ because $v$ has zero
trace there. To verify that $V\in H^1(D)$, choose a smooth lifting
$\eta$ of $f$ supported in a collar of the common outer arc away from
$F$. Extend $v-\eta|_U\in H^1_0(U)$ by zero, then add $\eta$.
The resulting $V$ and $\varphi$ have the same outer boundary trace.
This argument does not require regularity at the cusp of $F$.

Put $b_f(t)=\int_\alpha^\pi B_\kappa(t,s)f(s)\,\mathrm ds$.
The outward derivative of $v$ on $\Gamma$ is $-b_f$. Integration by
parts gives the transmission identity
$(-\Delta-\kappa^2)V=-b_f\mathcal H^1|_\Gamma$. Hence $\varphi-V$ has zero
boundary data and positive source $b_f\mathcal H^1|_\Gamma$ when
$f\geq0$. By the trace theorem on $U$, this source belongs to $H^{-1}(D)$.
Green representation, using $V=0$ on $F$, yields
\eqref{eq:eigen-transfer} by approximating a point mass at $s$.
The interface densities converge uniformly: the source is separated
from $\Gamma$, and odd reflection across the diameter gives
$B_\kappa(t,s)=O(t)$ at $t=0$. At the tangency endpoint the source remains
separated and ordinary local boundary estimates apply.

For $z\in F$ the Green weight in \eqref{eq:eigen-transfer} is strictly
positive for $0<t<\alpha$. Applying Lemma~\ref{lem:eigen-kernel} under the integral gives
\[
 0<P_{D,\kappa}(z,s')\leq\frac{\pi-s'}{\pi-s}P_{D,\kappa}(z,s)
   <P_{D,\kappa}(z,s).
\]
At a fixed point of the open smooth outer cap boundary the target is
a positive distance from $\Gamma$. Boundary estimates justify taking
its inward derivative under the integral. Hopf's lemma makes the
derivative of the Green weight strictly positive, proving the same
comparison for the derivative kernels.
\end{proof}

\begin{theorem}[Angular balancing for the first eigenvalue]
\label{thm:eigen-balance}\label{prop:straight-balance}
At each fixed radius in the respective ranges
\[
 \begin{aligned}
 \frac13<r<\sqrt2-1&\quad\text{for the straight family }\eqref{eq:straight-caps},\\
 r_0<r<\frac12&\quad\text{for the curved family }\eqref{eq:gaps},
 \end{aligned}
\]
the equilateral cap body uniquely maximizes the first Dirichlet eigenvalue,
up to rotations. Thus, in either family,
\begin{equation}\label{eq:eigen-balance}
 \lambda_1(S_r)\leq\lambda_1(E_r),
\end{equation}
with equality precisely for equally spaced axes.
\end{theorem}
\begin{proof}
Use the notation of Theorem~\ref{prop:balance}: $\gamma=\alpha$ or
$\beta$, the axes are $0,b,-c$, with $2\gamma\leq b<c<\pi$ and
$b+c+2\gamma\leq2\pi$. Let $S$ be either cap body, normalize its positive
first eigenfunction by $\int_S\varphi^2=1$, and set $\lambda=\lambda_1(S)$.
The geometric reflection across the line at angle $-(c-b)/2$ is the one
already used in the torsion proof. Its comparison domain omits a nonempty
open part of $S$, so its first Dirichlet eigenvalue exceeds $\lambda$.
The maximum principle for $-\Delta-\lambda$ therefore gives
\begin{equation}\label{eq:eigen-remote}
 \varphi(re^{i(b+t)})\geq\varphi(re^{-i(c+t)}),
 \qquad -\gamma<t<\gamma.
\end{equation}
As before, the whole source arc at $b$ lies on the positive side of the
reflection line because $b+c+2\gamma\leq2\pi$.

Let $D^+$ be the upper half-disk with the upper half of the cap at zero
attached, and put $w(z)=\varphi(z)-\varphi(\bar z)$ there.
It solves $(-\Delta-\lambda)w=0$ and vanishes on the real interval and
exterior cap branch. The proper inclusions $B_r\subsetneq S$ and
$D^+\subsetneq S$ each omit an open region and give
\[
 \lambda=\lambda_1(S)<\frac{j_{0,1}^2}{r^2},\qquad
 \lambda<\lambda_1(D^+).
\]
After rescaling $D^+$ by $r^{-1}$, set
\[
 \kappa=r\sqrt{\lambda_1(S)},\qquad \kappa^2=r^2\lambda_1(S).
\]
The first bound gives $\kappa<j_{0,1}$, and the second gives
$\kappa^2<\lambda_1(r^{-1}D^+)=r^2\lambda_1(D^+)$, exactly the two
spectral hypotheses of Lemma~\ref{lem:eigen-cap-transfer}.
Its geometric hypotheses hold for both half-cap attachments.
Write $P(z,s)$ for the resulting kernel on $D^+$ with angular source
measure. At a target in the upper half-cap it is positive and strictly
decreasing for $\gamma<s<\pi$.

Define
\[
 a(t)=\varphi(re^{i(b-\gamma+t)}),\qquad
 f(t)=\varphi(re^{-i(c-\gamma+t)}),\qquad 0<t<2\gamma.
\]
The first reflection gives $a\geq f>0$. If $c+\gamma\leq\pi$, both
source angles lie in $(\gamma,\pi)$ and
\[
 P(z,b-\gamma+t)a(t)\geq P(z,b-\gamma+t)f(t)
 >P(z,c-\gamma+t)f(t),\qquad 0<t<2\gamma.
\]
Integrating gives $w>0$ in the upper half-cap. In the straight-cap
crossing case $c+\alpha>\pi$, use exactly the source decomposition
\eqref{eq:straight-source-integral}, replacing the torsion traces by
those of $\varphi$ and $P_0$ by $P$. The first two terms are again
strictly ordered on $(0,m)$ and all remaining contributions are
nonnegative.
For fixed $z$ in the cap and $\gamma<s_0<\pi$,
Lemmas~\ref{lem:eigen-kernel} and~\ref{lem:eigen-cap-transfer} give
\[
 0<P(z,s)\leq\frac{\pi-s}{\pi-s_0}P(z,s_0),\qquad s_0<s<\pi.
\]
Consequently $P(z,s)\to0$ as $s\uparrow\pi$. This ensures absolute
integrability at the truncated endpoint for bounded source traces,
even when the separate traces have a jump and need not have finite energy.
Thus $w>0$ in the upper half-cap for both geometries, including the
crossing case and touching configurations.

Hopf's lemma gives $p_+(x)>p_-(x)>0$ on $\Gamma_+$, where
$p_+(x)=-\partial_n\varphi(x)$ and
$p_-(x)=-\partial_n\varphi(\sigma x)$.
The common paired velocities \eqref{eq:rotspeed} and the Dirichlet
Hadamard formula give
\begin{equation}\label{eq:eigen-angular-derivative}
 \frac{\partial\lambda_1(S)}{\partial a_0}
 =-\int_{\Gamma_+}V_n(x)(p_+^2-p_-^2)\dd s<0.
\end{equation}
A clockwise motion therefore increases the eigenvalue while increasing
the smaller gap and decreasing the larger one. It preserves all gap
bounds for sufficiently small positive time. The one-sided derivative
at contacts is justified in Appendices~\ref{app:variations}
and~\ref{app:straight}.

The admissible gap sets are compact. Hausdorff continuity, obtained
from the common inscribed disk, domain monotonicity and scaling, gives
an eigenvalue maximum. Every unequal gap triple admits the increasing
variation just proved, so precisely the equilateral configurations
maximize the first eigenvalue.
\end{proof}

\subsection{The spectral comparison on both hexagon families}
\label{sec:eigen-hexagon}

We now vary the cap parameter on either equilateral hexagon family. The geometry
and normal velocities are those of
Subsections~\ref{sec:hexagon-geometry}--\ref{sec:torsion-hexagon}; we must
compare eigenfunction fluxes at the paired points of each circular or straight edge.
Here no boundary-source kernel is needed. The reflection used for
torsion in Lemma~\ref{lem:flux} leads to a Helmholtz equation with
Dirichlet data on two sides and Neumann data on the third.
Dirichlet domain monotonicity does not by itself supply a spectral gap
for this mixed problem. We obtain that gap from Green's identity,
then apply the mixed maximum principle and the Hadamard formula.

\begin{proposition}\label{prop:eigen-hexagon}
The functions $r\mapsto\lambda_1(H_r^{\rm s})$ on $[1/3,1/2]$ and
$r\mapsto\lambda_1(H_r)$ on $[r_0,1/2]$ are strictly decreasing.
\end{proposition}
\begin{proof}
Let $H$ be either hexagon, with $r$ in the interior of its parameter
interval. Put $\lambda=\lambda_1(H)$ and normalize its positive first
eigenfunction by $\int_H\varphi^2=1$.
We use the sector, points $A,B,M,Q$, and reflection $\sigma$ introduced
in Lemma~\ref{lem:flux}, and set
\[
 Q=(\xi,0),\qquad \xi=\frac{1-2r}{2-3r},\qquad
 F=D_B,\quad \Gamma=BM,\quad\Sigma=QM,\quad N=QB.
\]
The first eigenfunction is log-concave \cite{bm_eigen}, so its superlevel
sets are convex. Simplicity gives the reflection symmetries of $H$. By Lemma~\ref{lem:gradient-cone},
\[
 \nabla\varphi=-c_0e_0-c_1e_{\pi/3},\qquad c_0,c_1\geq0.
\]
The outward normal on $N$ is
$\nu=(-\sqrt3r/2,r/2-\xi)/|QB|$. Consequently
\begin{equation}\label{eq:eigen-neumann-sign}
 \partial_\nu\varphi
 =\frac{\sqrt3(rc_0+\xi c_1)}{2|QB|}\geq0
 \qquad\text{on }N.
\end{equation}

The sign just obtained is the input to the spectral-gap argument.
Let $\mu$ be the first eigenvalue on $F$ with Dirichlet conditions on
$\Gamma\cup\Sigma$ and Neumann conditions on $N$, and let $z>0$ be
its first eigenfunction. We verify that the mixed corners are acute, so
that Green's identity applies. For the curved family, if
$\psi=\arg(B-Q)$, then
$r/2-\xi=-(3r^2-6r+2)/(2(2-3r))>0$ and $\xi>0$ imply
$\pi/3<\psi<\pi/2$. The angle at $Q$ is $\psi/2$.
If $\theta_B$ is the outward normal angle at $B$, the angle at $B$ is
$\pi/2+\theta_B-\psi<\pi/2$, because $\theta_B\leq\pi/3$.
For the straight family, $F=QBM$ is half of the isosceles triangle
$QAB$: the angle at $Q$ is half its apex angle and the angle at $B$
is a base angle, so both mixed corners are acute. In either family,
$M$ is a right-angle Dirichlet--Dirichlet corner.
The mixed-corner regularity in \cite[Chapter~4]{Grisvard} justifies the
normal derivatives and Green's identity, which gives
\begin{equation}\label{eq:mixed-spectral-gap}
 (\mu-\lambda)\int_F\varphi z
 =-\int_\Sigma\varphi\partial_\nu z
   +\int_N z\partial_\nu\varphi>0.
\end{equation}
Since $\varphi>0$ on the open median and $\partial_\nu z<0$ there by
Hopf's lemma, the first term is strictly positive. The second is
nonnegative by \eqref{eq:eigen-neumann-sign}. Since
$\int_F\varphi z>0$, we obtain $\mu>\lambda$. This proves the mixed
spectral gap before it is used in a maximum-principle argument.

The difference $w=\varphi-\varphi\circ\sigma$ satisfies
\[
 (-\Delta-\lambda)w=0\quad\text{in }F,\qquad
 w=0\quad\text{on }\Gamma\cup\Sigma,\qquad
 \partial_\nu w=\partial_\nu\varphi\geq0\quad\text{on }N.
\]
The last equality uses the Neumann symmetry condition on the reflected
segment $QA$. Testing with $w^- =\max\{-w,0\}$ and using the
nonnegative outward Neumann data yields
\[
 \int_F|\nabla w^-|^2-\lambda\int_F(w^-)^2\leq0.
\]
The mixed Rayleigh principle bounds the left side below by
$(\mu-\lambda)\|w^-\|_2^2$. Thus $w^-=0$, so $w\geq0$.
Comparison on the open arc or segment $BM$ gives
\begin{equation}\label{eq:eigen-paired-flux}
 p_B\geq p_A>0,\qquad p=-\partial_n\varphi,
\end{equation}
where the subscripts indicate paired points closer to the respective
endpoints of the boundary piece.

The flux order is now the same as in the torsion proof; only the sign
in the Hadamard formula changes. With the
circle data in \eqref{eq:center}, for the curved family we write
$X(\pm s)=C+\cos s\,n_m\pm\sin s\,t_m$, where
$0<s<\alpha/2$. The corresponding normal velocities are
$V_\pm=a_r\cos s\pm b_r\sin s$, with the positive coefficients
$a_r,b_r$ from \eqref{eq:velocity}; the plus sign is the $B$-half.
For this family, Hadamard's formula and the six symmetries give
\begin{equation}\label{eq:eigen-hexagon-derivative}
 \frac{\dd}{\dd r}\lambda_1(H_r)
 =-6\int_0^{\alpha/2}\left[
 a_r\cos s\,(p_B^2+p_A^2)
 +b_r\sin s\,(p_B^2-p_A^2)\right]\dd s<0.
\end{equation}
The first term in brackets is strictly positive; the second is
nonnegative by \eqref{eq:eigen-paired-flux}. Corollary~\ref{cor:app-hexagon-variation} justifies this variation and
gives continuity at both endpoints.
For the straight family, the affine normal speed $b$ on each edge has
$b(U)>0$ and $b(U)+b(V)>0$ at the paired points, as proved in
Proposition~\ref{prop:hexmono}. Hence
\[
 p_B^2b(U)+p_A^2b(V)
 =(p_B^2-p_A^2)b(U)+p_A^2\{b(U)+b(V)\}>0.
\]
The negative sign in the eigenvalue Hadamard formula gives
$\lambda_1'(H_r^{\rm s})<0$. The piecewise affine sector maps used
there justify this variation. Continuity at the endpoint polygons
completes both assertions on their closed intervals.
\end{proof}

\subsection{Proof of the eigenvalue width inequalities}
\label{sec:eigenvalue-final}\label{sec:straight-analysis}
Subsections~\ref{sec:eigen-angular} and~\ref{sec:eigen-hexagon} establish
the additional spectral ingredients needed to apply the geometry of
Section~\ref{sec:torsion} to the first eigenvalue. Table~\ref{tab:analytic-dictionary}
summarizes the two analytic arguments after those ingredients have been
proved; $\int_\Omega\varphi^2=1$ and $V_n=V\cdot n$ in the table.
\begin{table}[H]
\centering
\small
\renewcommand{\arraystretch}{1.12}
\begin{tabularx}{\linewidth}{@{}>{\raggedright\arraybackslash}X >{\raggedright\arraybackslash}X@{}}
\toprule
\textbf{Torsion} & \textbf{First eigenvalue}\\
\midrule
$-\Delta u=1$ & $-\Delta\varphi=\lambda_1(\Omega)\varphi$\\
\addlinespace[.4ex]
Reflected difference is harmonic
& Reflected difference solves $(-\Delta-\lambda_1(\Omega))w=0$\\
\addlinespace[.4ex]
Laplace maximum principle applies directly
& Requires $\lambda_1(\Omega)<\lambda_1(D)$, or a mixed spectral gap\\
\addlinespace[.4ex]
Explicit harmonic boundary kernel $B_0$
& Helmholtz boundary kernel $B_\kappa$\\
\addlinespace[.4ex]
Kernel monotonicity by explicit differentiation
& Kernel monotonicity from a positive separated Green expansion\\
\addlinespace[.4ex]
Gradient cone from symmetry and concavity of $\sqrt u$
& Gradient cone from symmetry and log-concavity of $\varphi$\\
\addlinespace[.4ex]
$DT(\Omega)[V]=\int_{\partial\Omega}(\partial_nu)^2V_n\dd s$
& $D\lambda_1(\Omega)[V]=-\int_{\partial\Omega}(\partial_n\varphi)^2V_n\dd s$\\
\addlinespace[.4ex]
Harmonic mixed comparison on the hexagon sector
& Mixed Helmholtz comparison after proving $\mu>\lambda_1(\Omega)$\\
\bottomrule
\end{tabularx}
\caption{Analytic dictionary between the torsion and first-eigenvalue
reflection arguments. The geometric reflections are the same; the
eigenvalue problem requires additional spectral-gap and Helmholtz-kernel inputs.}
\label{tab:analytic-dictionary}
\end{table}

\subsubsection{The inscribed-disk bound at minimal width}
\label{sec:straight-disk-bounds}
For a body $K$ of minimal width one and inradius $r$,
\[
 \lambda_1(K)\leq\frac{j_{0,1}^2}{r^2}.
\]
The eigenvalue cutoff is
\begin{equation}\label{eq:straight-cutoffs}
 r_\lambda=\frac{j_{0,1}}{2\pi}<\frac25<\sqrt2-1,
\end{equation}
using $j_{0,1}<5/2$ and $\pi>25/8$.
Thus $\lambda_1(K)<\lambda_1(\mathcal T)$ whenever $r>r_\lambda$.
At $r=r_\lambda$, strictness follows from proper inclusion of the
inscribed disk, whose width $2r_\lambda$ is less than one.
The omitted open region gives strict Dirichlet eigenvalue monotonicity,
as justified in Subsection~\ref{sec:eigen-angular}. It remains to treat
$1/3<r<r_\lambda$, within the straight-cap range of
Theorem~\ref{thm:eigen-balance}.

\subsubsection{Completion of the two width comparisons}
\label{sec:curved-eigenvalue-final}\label{sec:straight-final}
At minimal width, for $1/3<r<r_\lambda$, the straight-cap reduction,
Theorem~\ref{thm:eigen-balance}, the inclusion
$H_r^{\rm s}\subset E_r^{\rm s}$, and
Proposition~\ref{prop:eigen-hexagon} give
\[
 \lambda_1(K)\leq\lambda_1(S_r^{\rm s})\leq\lambda_1(E_r^{\rm s})
 \leq\lambda_1(H_r^{\rm s})<\lambda_1(H_{1/3}^{\rm s})
 =\lambda_1(\mathcal T).
\]
The disk estimate covers larger inradii, and at $r=1/3$
Lemma~\ref{lem:straight-reduction} gives $K=\mathcal T$.

At constant width, for $r_0<r<1/2$, domain monotonicity,
Theorem~\ref{thm:eigen-balance}, Lemma~\ref{lem:hexgeom}, and
Proposition~\ref{prop:eigen-hexagon} give
\[
 \lambda_1(K)\leq\lambda_1(S)\leq\lambda_1(E_r)
 \leq\lambda_1(H_r)<\lambda_1(H_{r_0})=\lambda_1(\mathcal R).
\]
If $r=1/2$, $K=B_{1/2}$ contains $H_{1/2}$ and the same strict bound
follows. At $r=r_0$, the geometric equality argument in
Subsection~\ref{sec:torsion-final} forces $K=\mathcal R$.
Finally, $\lambda_1(w\Omega)=w^{-2}\lambda_1(\Omega)$ proves both
spectral assertions and their equality cases in Theorem~\ref{thm:main}.

\section{Integral functionals and interaction kernels}
\label{sec:integral-kernels}
We use the cap subsets to study radial integrals and interaction energies.

\subsection{Radial integrals and interaction energies}\label{sec:kernels}
We consider planar convex bodies of minimal width one, or constant
width one. Let $Q$ denote the {\bf equilateral triangle of altitude one in
the first class} and {\bf the Reuleaux triangle of width one in the second}.
Their inradii are $1/3$ and $1-1/\sqrt3$, respectively; in either case put $b=r(Q)$.
All one-point comparisons below place the origin at an incenter of
the original body. At constant width this is the common center of
the smallest enclosing circle and largest inscribed circle. 

For a nonnegative radial kernel $k$, define
\[
 \J_k(K,z)=\int_K k(|x-z|)\dd x,\qquad
 \E_k(K)=\int_K\int_K k(|x-y|)\dd x\dd y.
\]
For a kernel $H$ on $\R^2$, we also write
$\E_H(K)=\iint_{K\times K}H(x-y)\dd x\dd y$; the radial notation
means $H(z)=k(|z|)$. The energies are unnormalized Lebesgue integrals. Dividing by $|K|$
or $|K|^2$ changes their inclusion properties and is a different problem.
For Riesz kernels we use $k_s(t)=t^{-s}$, $0<s<2$; both integrals
are then finite on bounded sets.

Bogosel, Bucur and Fragal\`a \cite{BogoselBucurFragala} prove a polygonal
Hardy--Littlewood inequality and investigate two-point energies at
prescribed area and number of sides. They also show that the kernel
can affect symmetry in the polygonal Riesz problem. Here the area
varies, and the relevant extremum is a minimum. Nevertheless, the
reduction to disk-indicator kernels and the comparison of simpler
subsets remain useful.

We apply the two geometric reductions of
Subsections~\ref{sec:straight-geometry} and~\ref{sec:geometry}. To treat
them together, we use $a$ for their opening half-angle: thus $a$ is
$\alpha$ in the straight case and $\beta$ in the curved case.

\subsubsection{The two cap families}
Put $e_\theta=(\cos\theta,\sin\theta)$.
The three cap vertices are $(1-r) e_{a_j}$. Write $a$ for the half-angle of
one cap, and let $\delta_j$ be the cyclic gaps between apex axes.
For straight caps,
\begin{equation}\label{kr:eq:straight}
 \begin{gathered}
 \tfrac13\leq r<\tfrac12,\qquad
 a=\arccos(r/(1-r)),\\
 S_r=\conv(\overline B_r,(1-r) e_{a_1},(1-r) e_{a_2},(1-r) e_{a_3}),
 \quad \delta_j\geq2a.
 \end{gathered}
\end{equation}
For circular caps,
\begin{equation}\label{kr:eq:curved}
 1-\tfrac1{\sqrt3}\leq r<\tfrac12,\qquad
 a=2\arccos\frac1{2(1-r)},\qquad
 2a\leq\delta_j\leq\pi-a.
\end{equation}
In this case $S_r$ is the union of $\overline B_r$ and the three
unit-ball hulls $\bh(\overline B_r\cup\{(1-r) e_{a_j}\})$.
These unions are convex and the cap interiors are disjoint.
In both cases $\sum_j\delta_j=2\pi$ and $0<a\leq\pi/3$.

After translation, every original body with inradius $r<1/2$
contains a corresponding $S_r$. The construction uses three contact
points of an inscribed circle in the straight case and three contact
points of the concentric enclosing circle in the circular case.
At $r=b$, the gap bounds force $\delta_j=2\pi/3$, and $S_b=Q$.
The equality cases in the radius bounds also force the original body
to be $Q$. At $r=1/2$ we use $\overline B_{1/2}\subset K$.
The containments were proved in Lemmas~\ref{lem:straight-reduction}
and~\ref{lem:geometry}; the endpoint rigidity in both classes is
proved in Theorem~\ref{thm:inner}.

For an axis at zero the radial boundary on $|\phi|\leq a$ is
\begin{equation}\label{kr:eq:radial}
 p_r(\phi)=
 \begin{cases}
 r\sec(a-|\phi|),&\text{straight},\\
 -(1-r)\cos(a-|\phi|)+\sqrt{1-(1-r)^2\sin^2(a-|\phi|)},&\text{circular}.
 \end{cases}
\end{equation}
It equals $r$ on the remaining circle sectors. In particular
$r\leq p_r\leq1-r$. Let $E_r$ denote the equally spaced configuration.
Every radial one-point integral on $S_r$ is independent of its cap
positions. For the two-point energy this is true only of the disk
term, disk--cap terms and individual cap self-interactions.

\subsubsection{Concentric clipping and the one-point inequality}
For a fixed radius $R\geq0$, put
\[
 C_r(R)=|S_r\cap B_R|.
\]
This quantity does not depend on the cap axes. We may therefore
differentiate the equilateral representative, which is admissible for
the whole interval of $r$. We prove monotonicity by examining the normal
velocity of its boundary. Recall the notation $b = r(Q)$, the lower bound on the inradius for each class: minimal width and constant width.

\begin{lemma}[One change of sign]\label{kr:lem:velocity}
For either cap family and $b<r<1/2$, the normal velocity under
increasing $r$ is positive at the smallest boundary radii and negative
at the largest boundary radii, with exactly one transition on each
half-cap. Its integral over the whole boundary is strictly positive.
\end{lemma}
\begin{proof}
For straight caps, take $0<\phi<a$ and set $z=a-\phi$. Differentiating
the first formula in \eqref{kr:eq:radial} at fixed polar angle gives
\[
 \partial_r p_r=\sec z\{1+r a'\tan z\},\qquad
 a'=-\frac1{(1-r)^2\sin a}.
\]
The sign of the normal speed is the sign of this derivative.
The expression in braces decreases strictly as $z$ goes from $0$
to $a$, and its endpoint values have signs $+$ and $-$.
The boundary radius $r\sec z$ increases strictly with $z$.
The remaining incircle arcs have speed one. The area is
\[
 A(r)=r^2\{\pi+3(\tan a-a)\},
\]
and direct differentiation gives
\begin{equation}\label{kr:eq:straight-area-derivative}
 A'(r)=r\{2\pi-6a+3\tan a(1-\cos a)\}>0.
\end{equation}

For circular caps, write the upper arc as
\[
 x(\theta)=-(1-r) e_a+e_\theta,\qquad a/2\leq\theta\leq a,
 \qquad a'=-\frac1{(1-r)^3\sin a}.
\]
Its outward normal speed is
\begin{equation}\label{kr:eq:curved-speed}
 V_n=\cos z-2\cot(a/2)\sin z,\qquad z=a-\theta\in[0,a/2].
\end{equation}
This is strictly decreasing in $z$. Its values at the join and apex
are $1$ and $-\cos(a/2)$, respectively. Moreover,
\[
 |x|^2=1+(1-r)^2-2(1-r)\cos z
\]
increases with $z$. The area and its derivative are
\begin{align}
 A(r)&=\pi r^2+\frac32(1-2r^2)a
                  -\frac32\sqrt{4(1-r)^2-1},\label{kr:eq:curved-area}\\
 A'(r)&=r(2\pi-6a)+6\{\sin(a/2)-2\cot(a/2)(1-\cos(a/2))\}>0.
 \label{kr:eq:curved-area-derivative}
\end{align}
For the last sign, the first term is nonnegative and
$\cot(a/4)>2\cot(a/2)$. These computations include all boundary
pieces; the moving endpoints contribute no separate term because
the incident pieces have the same endpoints.
\end{proof}

\begin{theorem}[Clipped-area comparison]\label{kr:thm:clipping}
For both cap families, $r\mapsto C_r(R)$ is nondecreasing on
$[b,1/2]$ for every $R\geq0$. It is strictly increasing wherever
$b<r<R$. Consequently, if $K$ is placed at an incenter,
\begin{equation}\label{kr:eq:clipping}
 |K\cap B_R|\geq |Q\cap B_R|\qquad(R\geq0).
\end{equation}
For $R\leq b$ the two sides equal $\pi R^2$ for every $K$.
For each fixed $R>b$, equality characterizes $Q$.
\end{theorem}
\begin{proof}
There is a radius $\rho_*\in(r,1-r)$ such that $V_n>0$ below
$\rho_*$ and $V_n<0$ above it. Except at a tangency value, shape
differentiation of an intersection with the fixed disk gives
\begin{equation}\label{kr:eq:clipping-derivative}
 \partial_r C_r(R)=\int_{\partial S_r\cap B_R} V_n\dd s.
\end{equation}
This also follows by differentiating the radial integral
$\frac12\int\min(p_r(\phi)^2,R^2)\dd\phi$.
If $R\leq r$, $C_r(R)=\pi R^2$. If $r<R\leq\rho_*$,
the integral in \eqref{kr:eq:clipping-derivative} is strictly positive.
If $\rho_*<R<1-r$, it is at least the full boundary integral $A'(r)>0$,
since only negative speeds have been omitted. For $R\geq1-r$ it is
exactly $A'(r)>0$. The functions are locally absolutely continuous
in $r$, so these signs and continuity give the assertion, including
the isolated tangencies and endpoint configurations.

Now use $S_r\subset K$ and compare with $S_b=Q$. When $R>b$ and
$r>b$, the interval $(b,\min(r,R))$ has positive length, giving a
strict inequality. The case $r=1/2$ follows by continuity and disk
inclusion. If equality holds, $r=b$ and the geometric rigidity gives
$K=Q$ up to rotation. For $R\leq b$, both bodies contain $B_R$.
\end{proof}

\begin{corollary}[Hardy--Littlewood-type minimum]\label{kr:cor:HL}
Let $k\geq0$ be nonincreasing and locally integrable as a radial
function on $\R^2$. For a planar body of minimal width one or constant
width one, and an incenter $O$,
\begin{equation}\label{kr:eq:HL}
 \J_k(K,O)\geq\J_k(Q,0).
\end{equation}
Write $k(t)=L+\int_{(t,\infty)}\dd\mu(R)$ with $L\geq0$ and
$\mu$ a nonnegative Stieltjes measure. If
$L>0$ or $\mu((b,\infty))>0$, equality characterizes $Q$.
In particular this applies to every one-point Riesz kernel
$k(t)=t^{-s}$, $0<s<2$.
\end{corollary}
\begin{proof}
Tonelli's theorem gives
\[
 \J_k(K,O)=L|K|+\int_{(0,\infty)}|K\cap B_R(O)|\dd\mu(R).
\]
Apply Theorem~\ref{kr:thm:clipping}, including its large-$R$ area case.
The stated condition puts positive weight on a strict comparison
unless $K=Q$. Values at discontinuity radii do not affect the integrals.
\end{proof}

The choice of an incenter is part of the assertion. A prescribed
Steiner point or centroid would require a separate argument.
An unrestricted translation minimization for a decreasing kernel is
also a different problem: a body can be moved far from the origin.

\begin{corollary}[Optimization over the center]\label{kr:cor:max-overlap}
For every $R\geq0$ and either width class,
\begin{equation}\label{kr:eq:max-overlap}
 \max_{z\in\R^2}|K\cap B_R(z)|
 \geq\max_{z\in\R^2}|Q\cap B_R(z)|=|Q\cap B_R(0)|.
\end{equation}
For $R\leq b$ the maximum equals $\pi R^2$ for every $K$;
for each $R>b$, equality in \eqref{kr:eq:max-overlap} characterizes $Q$.
More generally, for a kernel as in Corollary~\ref{kr:cor:HL},
\[
 \sup_z\J_k(K,z)\geq\sup_z\J_k(Q,z)=\J_k(Q,0).
\]
The same strictness condition as in Corollary~\ref{kr:cor:HL} gives
uniqueness of the minimizing shape.
\end{corollary}
\begin{proof}
The left side of \eqref{kr:eq:max-overlap} is at least the value at an
incenter, so Theorem~\ref{kr:thm:clipping} supplies its lower bound.
It remains to show that zero maximizes the overlap of $Q$.
Let $\omega$ denote rotation through $2\pi/3$ about zero.
The three sets $Q\cap B_R(\omega^j z)$, $j=0,1,2$, are rotations
of one another, and their Minkowski average is contained in
$Q\cap B_R(0)$. The Brunn--Minkowski inequality therefore gives
\[
 |Q\cap B_R(0)|^{1/2}\geq
 \frac13\sum_{j=0}^2|Q\cap B_R(\omega^j z)|^{1/2}
 =|Q\cap B_R(z)|^{1/2}.
\]
Zero-overlap sets need no separate estimate. This proves the equality
on the right and, by layer cake, the optimized-kernel statement.
The strictness assertion follows from the incenter lower bound.
\end{proof}

There is no corresponding statement for all increasing kernels.
For example $k(t)=\ind_{(1/2,\infty)}(t)$ has zero integral on
$B_{1/2}$, whereas neither reference triangle fits into a disk of
radius $1/2$, whatever its translation. Their integrals are positive.

\subsubsection{Angular balancing for convex decreasing pair kernels}
The inclusion $S_r\subset K$ gives $\E_k(K)\geq\E_k(S_r)$ for
every nonnegative kernel. We now give a second comparison for convex
decreasing kernels. Proposition~\ref{kr:prop:cutoff} below shows that
nonnegativity and monotonicity alone do not imply this angular comparison.

Let $A$ be the portion of the axis-zero cap outside $B_r$.
Write a point of it as $\rho e_\phi$, where $|\phi|\leq a$ and
$r\leq\rho\leq p_r(\phi)$. Both geometries satisfy
\begin{equation}\label{kr:eq:tangent-bounds}
 \rho\cos(a-|\phi|)\leq r,\qquad (1-r)/r\leq\sec a.
\end{equation}
The first inequality is the supporting-line inequality at the nearby
incircle tangency. The second is equality for straight caps; for
circular caps it follows from $(1-r)\cos a=1/(2(1-r))-(1-r)\leq r$.

\begin{lemma}[Distances between separate caps]\label{kr:lem:distance}
Fix $\rho e_\phi,\sigma e_\psi\in A$. For
$2a\leq\delta\leq2\pi-2a$, the function
\[
 D(\delta)=|\rho e_\phi-\sigma e_{\delta+\psi}|
\]
is concave wherever the points are distinct.
\end{lemma}
\begin{proof}
Let $\vartheta=\delta+\psi-\phi\in[0,2\pi]$, and put
$q=\max(\rho,\sigma)/\min(\rho,\sigma)$. We claim
\begin{equation}\label{kr:eq:cosine-condition}
 \cos\vartheta\leq1/q.
\end{equation}
It is immediate if $\cos\vartheta\leq0$. Otherwise let $t$ be the
smaller angular distance, $0\leq t<\pi/2$. If $t\geq a$, use
$q\leq (1-r)/r\leq\sec a\leq\sec t$.
If $t<a$ and $\vartheta=t$, put $u=a-\phi$, $v=a+\psi$.
Then $u,v\geq0$ and $u+v\leq t<a$. Thus the two tangent
inequalities give $\rho,\sigma\leq r\sec t$, whereas both are
at least $r$, proving $q\leq\sec t$ again. If
$\vartheta=2\pi-t$, use the other two tangencies in the same way.
This proves \eqref{kr:eq:cosine-condition}.

Writing $c=\cos\vartheta$, direct differentiation gives
\[
 D''=\frac{\rho\sigma}{D^3}
 \{(\rho^2+\sigma^2)c-\rho\sigma(1+c^2)\}\leq0.
\]
Indeed, when $c>0$ the opposite of the expression in braces, divided
by $\rho\sigma$, is $(1-qc)(1-c/q)\geq0$; for $c\leq0$ the
sign is immediate.
\end{proof}

\begin{theorem}[Angular energy comparison]\label{kr:thm:balancing}
Suppose $k$ is nonnegative, convex and nonincreasing on $(0,\infty)$,
and locally integrable as a radial kernel in the plane. For either
cap family at fixed $b\leq r<1/2$,
\begin{equation}\label{kr:eq:balancing}
 \E_k(S_r)\geq\E_k(E_r).
\end{equation}
If $k$ is $C^2$ and strictly decreasing, the equally spaced
configuration is the only minimizer up to rotation. This includes
all Riesz kernels $t^{-s}$, $0<s<2$.
\end{theorem}
\begin{proof}
Define the interaction of two cap regions by
\[
 F(\delta)=\int_A\int_A
 k(|\rho e_\phi-\sigma e_{\delta+\psi}|)\dd x\dd y.
\]
By Lemma~\ref{kr:lem:distance}, each distance is concave in $\delta$.
The composition with a convex nonincreasing $k$ is convex, so $F$
is convex on $[2a,2\pi-2a]$. Endpoint continuity follows from
local integrability and translation/rotation continuity of the
indicators; a uniform potential bound by the integral of $k$ over
a fixed disk controls the possible singularity at a touching point.
Equivalently, first work strictly inside the gap interval and pass
to its endpoints.

Decomposing $S_r$ into its disk and three caps gives
\begin{equation}\label{kr:eq:energy-decomposition}
 \E_k(S_r)=C(r)+2\{F(\delta_1)+F(\delta_2)+F(\delta_3)\},
\end{equation}
where $C(r)$ contains all self-interactions and disk--cap interactions.
The equality of the third pair interaction with $F(\delta_3)$ uses
$F(2\pi-\delta)=F(\delta)$. Since $\sum_j\delta_j=2\pi$,
Jensen's inequality proves \eqref{kr:eq:balancing}.
For a strictly decreasing $C^2$ kernel, $D''<0$ for almost every
point pair in the cap interiors. Thus $F$ is strictly convex on
the open interval, and equality forces all three gaps to be equal.
\end{proof}

\subsubsection{A decreasing non-convex cutoff kernel}
The convexity assumption in Theorem~\ref{kr:thm:balancing} cannot be
replaced by monotonicity alone. Consider
\[
 k(t)=\ind_{[0,0.95]}(t),
\]
which counts point pairs at distance at most $0.95$. This kernel is
nonnegative and nonincreasing, but its downward jump at $0.95$ makes
it non-convex. For this kernel, an unequal cap configuration can have
strictly smaller interaction energy than the equally spaced configuration.

\begin{proposition}\label{kr:prop:cutoff}
For each cap type there is a non-equilateral admissible $S_r$ with
\[
 \E_k(S_r)<\E_k(E_r).
\]
\end{proposition}
\begin{proof}
Set $r=499/1000$, so $1-r=501/1000$. For either type, its angle satisfies
$0<a<\pi/24$. Choose gaps
\[
 \delta_1=\delta_2=\frac\pi2+\frac a2,\qquad
 \delta_3=\pi-a.
\]
These satisfy all the straight and circular gap inequalities.
The three caps have the same positive area $m$.

For all equal-gap pairs and for the two shorter unequal-gap pairs,
the radial coordinates belong to $[r,1-r]$ and the smaller angular
separation is at most $2\pi/3+2a$. Therefore their point distances
are at most
\[
 \sqrt{(1-2r)^2+4(1-r)^2\sin^2(\pi/3+a)}
 <\sqrt{(1/500)^2+4(501/1000)^2\sin^2(3\pi/8)}
 <0.95.
\]
For the remaining pair, the smaller angular separation is at least
$\pi-3a$. All distances are at least
\[
 2r\cos(3a/2)>\frac{499}{500}\cos(\pi/16)>0.95.
\]
All inequalities are strict elementary inequalities; the last
trigonometric values can also be written in radicals.
Since $k$ equals one below the cutoff and zero above it, each of the
three cap-pair interactions in the equal-gap configuration is $m^2$.
In the unequal configuration two cap pairs still contribute $m^2$ each,
while the remaining pair lies entirely beyond $0.95$ and contributes zero.
The cap self-interactions and disk--cap interactions are unchanged.
Thus \eqref{kr:eq:energy-decomposition} gives
\[
 \E_k(S_r)-\E_k(E_r)=2(2m^2-3m^2)=-2m^2<0.
\]
\end{proof}

The strict distance margins also give the same energy ordering for
smooth nonincreasing approximations of this cutoff. Thus smoothness
alone does not replace convexity in Theorem~\ref{kr:thm:balancing}.
\subsection{Nonradial perturbations}\label{sec:nonradial}
Let $Q$ be the reference body, $b=r(Q)$ and $S_r$ a cap subset in the
corresponding width class. The area excess of $S_r$ controls its symmetric
difference from a suitably oriented copy of $Q$.

For a constant kernel $c>0$,
\[
 \int_S c=c|S|,\qquad \E_c(S)=c|S|^2,
\]
so the area theorem gives both minimization results. For a perturbed
kernel $c+q$, the constant part gives a gain proportional to
$|S_r|-|Q|$, while the perturbation error is controlled by
$|S_r\mathbin\triangle Q_\theta|$. Lemma~\ref{kr:lem:stability}
bounds this symmetric difference linearly by the same area excess,
allowing sufficiently small perturbations to be absorbed by the
constant-kernel gain, even without rotational symmetry.

\begin{lemma}[Stability inside the cap family]\label{kr:lem:stability}
For either family there is a finite constant $C>0$ such that every
admissible cap body $S_r$ admits a rotation $Q_\theta$ of $Q$ about
the origin satisfying
\begin{equation}\label{kr:eq:stability}
 |S_r\mathbin\triangle Q_\theta|\leq C\{|S_r|-|Q|\}.
\end{equation}
This linear estimate is uniform both in $r$ and in the admissible
angular positions of the three caps.
\end{lemma}
\begin{proof}
Let $\varepsilon=\pi/3-a$. Admissibility gives
$\delta_j\geq2a=2\pi/3-2\varepsilon$. Since the other two gaps
are also at least $2a$ and the three gaps sum to $2\pi$,
$\delta_j\leq2\pi-4a=2\pi/3+4\varepsilon$. Thus
\[
 -2\varepsilon\leq\delta_j-2\pi/3\leq4\varepsilon.
\]
In coordinates where the first cap axis is zero, the second and third
axes are $\delta_1$ and $2\pi-\delta_3$. Hence they differ from
$2\pi/3$ and $4\pi/3$, respectively, by at most $4\varepsilon$.
Let $Q_\theta$ denote the reference body with this alignment in the
original coordinates. Both angle functions satisfy $a(b)=\pi/3$
and are smooth with finite derivative near $b$, so
\[
 \varepsilon=\pi/3-a(r)=O(r-b).
\]

Near this endpoint, the cosine denominator in the straight profile
\eqref{kr:eq:radial} is at least $1/2$, and the square-root argument
in the circular profile stays uniformly bounded away from zero.
Thus each individual cap profile, extended by $r$ outside its sector,
depends uniformly Lipschitz-continuously on $r$, its opening angle
and its axis. Composing with absolute value and taking the maximum
of the three profiles preserve these uniform Lipschitz bounds.
The preceding parameter estimates therefore give
\[
 \|p_{S_r}-p_{Q_\theta}\|_\infty\leq C_0(r-b),
 \qquad d_H(S_r,Q_\theta)\leq C_1(r-b),
\]
where the Hausdorff estimate follows by comparison along each ray.

All these bodies have uniformly bounded perimeter. If convex planar
sets satisfy $d_H(K,L)\leq h$, then $K\subset L+hB_1$ and
$L\subset K+hB_1$; Steiner's formula gives
$|K\setminus L|\leq hP(L)+\pi h^2$ and
$|L\setminus K|\leq hP(K)+\pi h^2$.
Adding these inequalities yields, for $h=d_H(S_r,Q_\theta)$,
\[
 |S_r\mathbin\triangle Q_\theta|
 \leq\{P(S_r)+P(Q)\}h+2\pi h^2.
\]
Consequently $|S_r\mathbin\triangle Q_\theta|\leq C_2(r-b)$
near $b$, uniformly in the cap positions.

Write $A(r)=|S_r|$, which is independent of the cap positions.
The endpoint area derivatives are
\[
 A'_{\rm s}(1/3)=\frac{\sqrt3}{2}>0,\qquad
 A'_{\rm c}(1-1/\sqrt3)=21-12\sqrt3>0.
\]
By continuity, there is $c_1>0$ such that, for $r$ close to $b$,
\[
 A(r)-A(b)=\int_b^r A'(t)\dd t\geq c_1(r-b).
\]
Combining the two first-order estimates gives
\[
 |S_r\mathbin\triangle Q_\theta|
 \leq\frac{C_2}{c_1}\{A(r)-A(b)\},
\]
which proves \eqref{kr:eq:stability} near $b$.

Fix a sufficiently small $\eta>0$. On the remaining range
$r\in[b+\eta,1/2]$, strict area monotonicity gives
\[
 A(r)-A(b)\geq A(b+\eta)-A(b)>0,
\]
while the symmetric differences remain uniformly bounded. Hence the
ratio $|S_r\mathbin\triangle Q_\theta|/(A(r)-A(b))$ is uniformly
bounded on this compact parameter range as well. Together with the
estimate near $b$, this proves the lemma, including the disk endpoint
$r=1/2$ by continuity.
\end{proof}

\begin{remark}[Size of the stability constant]\label{kr:rem:stability-size}
The proof above is intended only to establish a uniform finite constant
$C$, and tracking its constants through the Hausdorff and Steiner
estimates would give a rather nonoptimal value. Since the cap profiles
are explicit, a more informative constant can instead be studied directly
through
\[
 C_* = \sup_{r,\,(a_1,a_2,a_3)}
 \frac{\displaystyle\inf_\theta|S_r\mathbin\triangle Q_\theta|}
      {\displaystyle|S_r|-|Q|},
\]
where the supremum is taken over the admissible cap configurations in
the corresponding family with $b<r\leq1/2$. The denominator depends
only on $r$, while the numerator can be evaluated from the radial
profiles in \eqref{kr:eq:radial}.

Let $C_{\rm align}$ denote the supremum of the same quotient when
$Q_\theta$ is chosen by aligning the first cap axis with a reference
axis, as in the proof, rather than minimizing over rotations.
Then $C_*\leq C_{\rm align}$. A direct numerical exploration with this
alignment suggests values of order
\[
 C_{\rm align}^{\rm s}\approx6.8,\qquad
 C_{\rm align}^{\rm c}\approx16.3
\]
for the straight and circular cap families, respectively, with the
largest ratios apparently occurring near the endpoint $r=b$.
These values concern the prescribed alignment and do not evaluate
the infimum over rotations in the definition of $C_*$.
Thus constants of the rough sizes $7$ and $17$ should be expected,
although we do not prove such numerical bounds here.

These values also indicate the scale of the perturbative conclusions
below. For example, if $C=7$ in the straight case and $C=17$ in the
circular case were valid stability constants, inserting them into
\eqref{kr:eq:exponent} would give sufficient Riesz ranges of approximately
\[
 s<2.9\times10^{-2},\qquad s<2.0\times10^{-2},
\]
respectively. These numbers are only indicative; no certified numerical
enclosure for $C_*$ is needed for the qualitative results of this section.
\end{remark}

Take $D=4/3$ for straight caps and $D=2/\sqrt3$ for circular caps.
All cap bodies and their reference bodies have diameter at most $D$.

For a nonradial kernel, we keep both the cap body and the kernel fixed
in the prescribed coordinates and choose the aligned comparison body
$Q_\theta$ supplied by Lemma~\ref{kr:lem:stability}. Only the
comparison body's orientation is chosen: all integrals use the same
fixed kernel in the original coordinates.

\begin{theorem}[Nonradial perturbations]\label{ck:thm:nonradial}
The following comparisons hold.
\begin{enumerate}
 \item Let $F\geq0$ be measurable on $\R^2$. On $B_{1-b}$ write $F=c+q$,
 with $c>0$ and
 $\|q\|_{L^\infty(B_{1-b}\setminus B_b)}<c/C$.
 Assume $F$ is integrable on $B_b$. If the original body is centered
 at an incenter, then
 \[
 \int_K F\geq\min_\theta\int_{Q_\theta}F.
 \]
 Every minimizing shape is a rotated reference body.
 \item Let $H(z)\geq0$ be a measurable interaction kernel, and write $H=c+q$
 on $B_D$, where $c>0$ and $D=2(1-b)$. If
 $\|q\|_{L^1(B_D)}<c|Q|/C$, then
 \[
 \iint_{K\times K}H(x-y)\dd x\dd y
 \geq\min_\theta\iint_{Q_\theta\times Q_\theta}H(x-y)\dd x\dd y.
 \]
 Every minimizing shape is a rotated reference body.
\end{enumerate}
\end{theorem}
\begin{proof}
Put $A=|S|$, $A_0=|Q|$ and choose $Q_\theta$ from
\eqref{kr:eq:stability}; we orient the reference body rather than
rotating the given kernel. For the first assertion, $S$ and $Q_\theta$
both contain $B_b$, so their difference avoids that disk. Hence
\[
 \int_S F-\int_{Q_\theta}F
 \geq(c-C\|q\|_\infty)(A-A_0)>0
\]
when $r>b$. The values on the common disk cancel even if singular.
For the second assertion only the even part of $H$ contributes to
the energy. Symmetrizing $q$ does not increase its $L^1$ norm.
Extend it by zero outside $B_D$. Its potentials against indicator
functions are bounded by $\eta=\|q\|_1$, and expansion gives
\[
 |\E_q(S)-\E_q(Q_\theta)|
 \leq2\eta|S\mathbin\triangle Q_\theta|
 \leq2C\eta(A-A_0).
\]
Adding the constant-kernel difference yields
$\{c(A+A_0)-2C\eta\}(A-A_0)>0$ for $r>b$.
Nonnegative inclusion transfers both comparisons to $K$.
The geometric endpoint $r=b$ forces the reference body; the
disk endpoint follows by continuity. Minimization over rotations
is attained by continuity of the integrals under rotation.
The norm condition concerns only the bounded cap subsets; nonnegative
inclusion applies even when the original minimal-width body has
arbitrarily large diameter.
\end{proof}

\begin{corollary}[A nonempty interval of Riesz exponents]\label{kr:cor:small-s}
For each width class there is $s_*>0$ such that $Q$ uniquely minimizes
$\E_s(K)=\iint_{K\times K}|x-y|^{-s}\dd x\dd y$ whenever
$0<s<s_*$. With $C,D$ as above, any
\begin{equation}\label{kr:eq:exponent}
 s_*\leq\frac{2|Q|}{\pi C D^2+|Q|}
\end{equation}
suffices. The present proof establishes finiteness of $C$ but does
not supply its numerical value.
\end{corollary}
\begin{proof}
Use $c=D^{-s}$ and $q(t)=(t^{-s}-D^{-s})\ind_{(0,D)}(t)$.
Then
\[
 \|q\|_1=2\pi\int_0^D(t^{-s}-D^{-s})t\dd t
          =\frac{\pi sD^{2-s}}{2-s}.
\]
Substituting into the second condition of Theorem~\ref{ck:thm:nonradial} gives exactly
$s<2|Q|/(\pi C D^2+|Q|)$.
\end{proof}

\appendix
\section{First variations at curved-cap contacts}
\label{app:variations}

We prove the first-variation formulas used for the two Dirichlet
functionals. We include admissible rotations which separate two
initially tangent caps, and variations of the circular hexagons.
For the cap rotation, we construct a map which follows the moving apex
in angular coordinates. Its derivative matrices have bounded first
difference quotients converging almost everywhere. These properties
suffice to differentiate the energies after pulling the problem back
to a fixed domain; differentiability in the $W^{1,\infty}$ norm is not
needed. We then recover the boundary formulas of \cite{HP} by integration
by parts.

Set $r_0=1-1/\sqrt3$, fix $r_0<r<1/2$, and write
\[
 \beta=\arccos\left(\frac{1}{2(1-r)^2}-1\right),
 \qquad e_\theta=(\cos\theta,\sin\theta),
 \qquad J(x_1,x_2)=(-x_2,x_1).
\]
For the cap with axis zero, the open upper branch and its normal are
\begin{equation}\label{eq:eigen-curved-cap}
 x(\theta)=-(1-r)e^{i\beta}+e^{i\theta},\qquad
 n=e^{i\theta},\qquad \beta/2<\theta<\beta.
\end{equation}
Thus $Jx\cdot n=(1-r)\sin(\beta-\theta)>0$, with the negative speed
on the reflected lower branch, as in \eqref{eq:rotspeed}.

The radial function of the disk together with the cap at $(1-r)e^{ia}$ is
\begin{equation}\label{eq:app-single-profile}
 p_a(\theta)=
 \begin{cases}
 -(1-r)\cos(\beta-t)+\sqrt{1-(1-r)^2\sin^2(\beta-t)},
      &t=\operatorname{dist}_{\mathbb S^1}(\theta,a)\leq\beta,\\
 r,&t\geq\beta.
 \end{cases}
\end{equation}
Here $(1-r)e^{ia}$ denotes $(1-r) e_a$ in complex notation.
Let $\Omega_0=S(r;a_0,a_1,a_2)$ be an admissible three-cap body. Choose
$\sigma\in\{-1,1\}$ such that rotating the first cap from $a_0$ to
$a_0+\sigma\tau$ is admissible for $0\leq\tau<\varepsilon$. Define
\[
 \Omega_\tau=S(r;a_0+\sigma\tau,a_1,a_2),\qquad
 p_\tau(\theta)=\max\{p_{a_0+\sigma\tau}(\theta),
                         p_{a_1}(\theta),p_{a_2}(\theta)\}.
\]
Thus $\sigma=1$ means counterclockwise rotation. Each domain contains the disk. Since the cap sectors have disjoint
interiors, each open exterior cap arc is a piece of its radial boundary
$p_\tau$.

\begin{lemma}[A deformation following the apex]\label{lem:cap-domain-map}
There are uniformly bi-Lipschitz maps
$F_\tau:\Omega_0\to\Omega_\tau$ and a Lipschitz vector field $V$ on
$\overline\Omega_0$ such that, as $\tau\downarrow0$,
\begin{align}\label{eq:app-map-limits}
 F_\tau&=I+\tau V+o(\tau)\quad\hbox{uniformly},\notag\\
 \|DF_\tau-I\|_{L^\infty}&\leq C\tau,
 &\frac{DF_\tau-I}{\tau}&\longrightarrow DV
                   \quad\hbox{almost everywhere}.
\end{align}
Their boundary velocity has normal component
\begin{equation}\label{eq:app-normal-speed}
 b=V\cdot n=
 \begin{cases}
 \sigma Jx\cdot n,&\text{on the two exterior arcs of the moving cap},\\
 0,&\text{on the remaining boundary},
 \end{cases}
\end{equation}
where the formula is understood at almost every boundary point.
The assertions remain valid when two caps are initially tangent.
\end{lemma}

\begin{proof}
We choose a smooth periodic function $\eta$ equal to one near $a_0$ and
zero near $a_1,a_2$. Choose these three neighborhoods disjoint and small
enough that the corresponding cap profile remains the uniquely active
one there for all sufficiently small admissible $\tau$.
The angular map
$\Theta_\tau(\theta)=\theta+\sigma\tau\eta(\theta)$ is an increasing
circle diffeomorphism when $\tau$ is small. Set
\begin{equation}\label{eq:app-radial-map}
 k_\tau(\theta)=\frac{p_\tau(\Theta_\tau(\theta))}{p_0(\theta)},
 \qquad
 F_\tau(\rho e_\theta)
   =\rho k_\tau(\theta)e_{\Theta_\tau(\theta)},
 \qquad F_\tau(0)=0.
\end{equation}
This map sends each radial interval bijectively onto the corresponding
radial interval of $\Omega_\tau$. Its denominator is bounded below by $r$.

Near the moving apex,
$p_\tau(\theta+\sigma\tau)=p_0(\theta)$; near either fixed apex,
$p_\tau(\theta)=p_0(\theta)$. Thus the change of angular variable in \eqref{eq:app-radial-map}
removes the moving jump in the derivative of the radial profile.
Away from these apex neighborhoods the profile pieces are smooth, with
uniformly bounded second angular derivatives up to their tangencies.
At a tangency angle $\theta_0$, its defining circle equation is
\[
 p^2+2(1-r)p\cos(\theta-\theta_0)+(1-r)^2=1.
\]
It gives, on the cap side,
\begin{equation}\label{eq:app-tangency-expansion}
 p(\theta)=r+\frac{r(1-r)}{2}(\theta-\theta_0)^2
                       +O(|\theta-\theta_0|^4),
 \qquad p'(\theta_0)=0.
\end{equation}
The value and first derivative therefore match those of the disk.
At a cap contact the two incident cap expressions also have the same
value and first derivative. If the contact opens, the interval where
the active boundary piece changes has length $O(\tau)$. On that
interval the radial remainder beyond its first-order displacement is
$O(\tau^2)$ and the angular derivative remainder is $O(\tau)$,
by \eqref{eq:app-tangency-expansion}. On all other pieces the formulas
depend smoothly on $\tau$. Combining the estimates on these finitely many pieces, we obtain
\begin{equation}\label{eq:app-k-bound}
 \|k_\tau-1\|_{W^{1,\infty}(\mathbb S^1)}\leq C\tau.
\end{equation}
In particular $F_\tau$ and its inverse are uniformly Lipschitz: the
angular maps and their inverses have bounded derivatives, and their
positive radial factors have bounded $W^{1,\infty}$ norms.

We next identify the first-order limit. Define $\dot p=-\sigma p'_0$ on the open
angular span of the moving cap and $\dot p=0$ elsewhere, disregarding
the finitely many junction angles. At either tangency $\dot p$ extends
continuously by zero. Put
\begin{equation}\label{eq:app-velocity-field}
 \begin{split}
 V(\rho e_\theta)&=\rho\left[
 \frac{\dot p(\theta)+\sigma\eta(\theta)p'_0(\theta)}{p_0(\theta)}
          e_\theta
       +\sigma\eta(\theta)J e_\theta\right],\\
 V(0)&=0.
 \end{split}
\end{equation}
The jump at the moving apex cancels in
$\dot p+\sigma\eta p'_0$, while this expression vanishes near the
fixed apex vertices. It is continuous at all tangencies and piecewise smooth
with bounded derivative. Thus the angular coefficients in
\eqref{eq:app-velocity-field} are Lipschitz, and their homogeneous
extension is Lipschitz also at the origin.
Using these piecewise formulas and the quadratic estimate
\eqref{eq:app-tangency-expansion}, we obtain the uniform expansion in
\eqref{eq:app-map-limits}. Off the finitely many initial junction rays,
the active profile is fixed for sufficiently small $\tau$ and is smooth
in both variables; hence the derivative quotients converge to $DV$.
The bound in \eqref{eq:app-map-limits} follows from
\eqref{eq:app-k-bound}. The exceptional rays and the origin have area zero. We therefore have
almost-everywhere convergence of the derivative quotients; their uniform
convergence is not needed.

To compute the normal boundary speed, we observe that a tangent vector
at a smooth radial boundary point $x=p_0(\theta)e_\theta$ is
$p'_0e_\theta+p_0J e_\theta$. On the moving cap,
\[
 V(x)-\sigma Jx
  =\sigma(\eta-1)(p'_0e_\theta+p_0J e_\theta),
\]
which is tangent. On the fixed boundary,
$V(x)=\sigma\eta(p'_0e_\theta+p_0J e_\theta)$ is tangent.
This proves \eqref{eq:app-normal-speed}, independently of the choice
of the interpolating function $\eta$.
\end{proof}

\begin{proposition}[Differentiating the pulled-back energies]
\label{prop:app-volume-derivatives}
Let $F_\tau$ be uniformly bi-Lipschitz maps satisfying
\eqref{eq:app-map-limits} on a bounded connected
Lipschitz domain $\Omega$, with $F_\tau(\Omega)=\Omega_\tau$.
Let $u$ be the torsion function on $\Omega$, and let $\varphi>0$ be its
first Dirichlet eigenfunction with
$\int_\Omega\varphi^2=1$, of eigenvalue $\lambda$.
The right derivatives exist and equal
\begin{align}\label{eq:app-torsion-volume}
 T'(0+)&=\int_\Omega
  \bigl[(2u-|\nabla u|^2)\operatorname{div}V
             +2\nabla u\cdot DV\nabla u\bigr],\\
 \lambda_1'(0+)&=\int_\Omega
  \bigl[(|\nabla\varphi|^2-\lambda\varphi^2)
                 \operatorname{div}V
            -2\nabla\varphi\cdot DV\nabla\varphi\bigr].
 \label{eq:app-eigen-volume}
\end{align}
\end{proposition}

\begin{proof}
We pull all functions back to $\Omega$. Since $DF_\tau$ is close to the
identity, its determinant is positive for small $\tau$. The coefficients
in the energy and mass integrals become
\[
 m_\tau=\det DF_\tau,
 \qquad A_\tau=m_\tau DF_\tau^{-1}DF_\tau^{-T}.
\]
They obey $\|A_\tau-I\|_\infty+\|m_\tau-1\|_\infty=O(\tau)$.
The matrices $A_\tau$ are uniformly elliptic, and $m_\tau$ is bounded
above and away from zero. Their first difference quotients are uniformly
bounded and converge almost everywhere to
\begin{equation}\label{eq:app-coefficient-limits}
 \dot m=\operatorname{div}V,
 \qquad \dot A=(\operatorname{div}V)I-DV-DV^T.
\end{equation}
Let
\[
 a_\tau(v,z)=\int_\Omega A_\tau\nabla v\cdot\nabla z,
 \qquad M_\tau(v,z)=\int_\Omega m_\tau vz,
 \qquad \ell_\tau(v)=\int_\Omega m_\tau v.
\]
By uniform ellipticity and Poincar\'e's inequality on the fixed domain,
$a_\tau$ is uniformly coercive on $H^1_0(\Omega)$.

For torsion, the pulled-back solution $U_\tau$ is the unique minimizer of
\[
 \mathcal J_\tau(v)=a_\tau(v,v)-2\ell_\tau(v),
 \qquad \min\mathcal J_\tau=-T(\Omega_\tau),
\]
and satisfies $a_\tau(U_\tau,z)=\ell_\tau(z)$. Subtract the equation
for $u=U_0$ and test with $U_\tau-u$. The coefficient bounds and
coercivity give
\[
 \|U_\tau-u\|_{H^1_0}\leq C\tau.
\]
Comparing each minimizer with the other one gives
\begin{align*}
 \mathcal J_\tau(U_\tau)-\mathcal J_0(U_\tau)
 &\leq\mathcal J_\tau(U_\tau)-\mathcal J_0(u)\\
 &\leq\mathcal J_\tau(u)-\mathcal J_0(u).
\end{align*}
After division by $\tau$, both bounds converge to
$\int_\Omega\dot A\nabla u\cdot\nabla u-2\dot m\,u$.
Indeed, the coefficient quotients are bounded and converge almost
everywhere; the strong $H^1$ convergence implies $L^1$ convergence of
the quadratic gradient densities. The negative of this common limit
is $T'(0+)$, giving \eqref{eq:app-torsion-volume}.

For the eigenvalue, write
$R_\tau(v)=a_\tau(v,v)/M_\tau(v,v)$ and let $\Phi_\tau>0$ be its
minimizer, normalized by $M_\tau(\Phi_\tau,\Phi_\tau)=1$.
Uniform comparison of the coefficients with $A_0=I,m_0=1$ gives
$\lambda_\tau\to\lambda$ and bounds $\Phi_\tau$ in $H^1_0$.
By weak compactness in $H^1_0$ and strong compactness in $L^2$, every
subsequential limit has $L^2$ norm one and minimizes $R_0$. Positivity
and simplicity of the first eigenvalue imply that this limit is
$\varphi$. Moreover
\[
 \int_\Omega|\nabla\Phi_\tau|^2
   =a_\tau(\Phi_\tau,\Phi_\tau)+O(\tau)
   =\lambda_\tau+O(\tau)\longrightarrow\lambda.
\]
The convergence of the gradient norms therefore gives strong
$H^1_0$ convergence of the entire family.

Comparing the Rayleigh quotients at the two minimizers, we obtain
\[
 R_\tau(\Phi_\tau)-R_0(\Phi_\tau)
 \leq\lambda_\tau-\lambda
 \leq R_\tau(\varphi)-R_0(\varphi).
\]
For either $v_\tau=\Phi_\tau$ or $v_\tau=\varphi$,
\[
 \frac{R_\tau(v_\tau)-R_0(v_\tau)}{\tau}
 =\frac{\displaystyle
 \frac{a_\tau-a_0}{\tau}(v_\tau,v_\tau)
 -R_0(v_\tau)\frac{M_\tau-M_0}{\tau}(v_\tau,v_\tau)}
 {M_\tau(v_\tau,v_\tau)}.
\]
By strong $H^1$ convergence and \eqref{eq:app-coefficient-limits}, both
bounds tend to
$\int_\Omega\dot A\nabla\varphi\cdot\nabla\varphi
                       -\lambda\dot m\,\varphi^2$.
This proves \eqref{eq:app-eigen-volume}. Since all comparisons were made for positive $\tau$, this proves the
assertions for the right derivatives.
\end{proof}

\begin{theorem}[Boundary formulas, including touching caps]
\label{thm:cap-first-variation}
For the admissible curved-cap rotation above, let $u$ be the torsion
function of $\Omega_0$ and normalize its positive first eigenfunction
$\varphi$ by $\int_{\Omega_0}\varphi^2=1$. Then
\begin{equation}\label{eq:app-boundary-derivatives}
 \left.\frac{\mathrm d}{\mathrm d\tau}T(\Omega_\tau)\right|_{0+}
   =\int_{\partial\Omega_0}(\partial_nu)^2b\,\mathrm ds,
 \qquad
 \left.\frac{\mathrm d}{\mathrm d\tau}\lambda_1(\Omega_\tau)
       \right|_{0+}
   =-\int_{\partial\Omega_0}(\partial_n\varphi)^2b\,\mathrm ds,
\end{equation}
where $b$ is given in \eqref{eq:app-normal-speed}. In coordinates in
which the moving cap has axis zero, its upper and lower arc speeds are
\[
 b_+(\theta)=\sigma(1-r)\sin(\beta-\theta),
 \qquad b_-(\theta)=-\sigma(1-r)\sin(\beta-\theta),
 \qquad \beta/2<\theta<\beta.
\]
There is no additional contribution from an apex or a disk--cap join.
\end{theorem}

\begin{proof}
We start from the volume formulas given by
Lemma~\ref{lem:cap-domain-map} and
Proposition~\ref{prop:app-volume-derivatives}. The domain is convex and piecewise smooth; hence the two
solutions belong to $H^2(\Omega_0)$ by convex-domain elliptic
regularity \cite[Chapter~3]{Grisvard}. Define the vector fields
\[
 \begin{split}
 Z_T&=(2u-|\nabla u|^2)V+2(V\cdot\nabla u)\nabla u,\\
 Z_\lambda&=(|\nabla\varphi|^2-\lambda\varphi^2)V
                         -2(V\cdot\nabla\varphi)\nabla\varphi.
 \end{split}
\]
They belong to $W^{1,1}(\Omega_0;\mathbb R^2)$: products of a Hessian
and a gradient are integrable by Cauchy--Schwarz, and $V$ is Lipschitz.
Using $-\Delta u=1$ and $-\Delta\varphi=\lambda\varphi$, direct
weak differentiation gives
\[
 \begin{split}
 \operatorname{div}Z_T
   &=(2u-|\nabla u|^2)\operatorname{div}V
                                  +2\nabla u\cdot DV\nabla u,\\
 \operatorname{div}Z_\lambda
   &=(|\nabla\varphi|^2-\lambda\varphi^2)\operatorname{div}V
                         -2\nabla\varphi\cdot DV\nabla\varphi.
 \end{split}
\]
Since the domain is Lipschitz, we can apply the divergence theorem to
these Sobolev vector fields. On its smooth boundary pieces, the zero Dirichlet
trace gives
$\nabla u=(\partial_nu)n$ and
$\nabla\varphi=(\partial_n\varphi)n$ almost everywhere. Their
normal traces are square integrable, by $H^2$ regularity. Therefore
\[
 Z_T\cdot n=(\partial_nu)^2b,
 \qquad Z_\lambda\cdot n=-(\partial_n\varphi)^2b.
\]
This proves \eqref{eq:app-boundary-derivatives}. For the simple-eigenvalue
formula see also \cite[Section~5.7]{HP}. Because the stress fields belong to $W^{1,1}$, the divergence theorem
accounts for the whole variation. The boundary integrals have
integrable arclength densities. Consequently there is no additional
atomic term at a corner or a disk--cap join.

Finally, on the upper arc $x=-(1-r) e_\beta+e_\theta$ and $n=e_\theta$,
whence $Jx\cdot n=(1-r)\sin(\beta-\theta)$. Reflection in the cap axis
reverses this scalar product and gives the stated lower speed.
\end{proof}

\begin{corollary}[The circular-hexagon parameter and endpoint continuity]
\label{cor:app-hexagon-variation}
For $r_0<r<1/2$, the family $H_r$ satisfies the boundary formulas
\eqref{eq:app-boundary-derivatives} with parameter $r$ and normal speed
$C'(r)\cdot n$ on an arc whose unit-circle center is $C(r)$.
Both functionals are continuous on $[r_0,1/2]$.
\end{corollary}

\begin{proof}
We use radial maps, since the six vertex rays are fixed. On each
sector between consecutive rays the radial boundary is the positive root
\[
 p_r(\theta)=C(r)\cdot e_\theta+
 \sqrt{1-|C(r)|^2+(C(r)\cdot e_\theta)^2}.
\]
Every defining disk contains the origin in its interior. On compact
subintervals of $(r_0,1/2)$ the root is uniformly regular, and the
profile is smooth in $r$ and $\theta$ up to each sector endpoint.
At a common endpoint the two profiles give the same moving vertex,
so their parameter derivatives agree there. With the angular map
fixed, the radial maps
\[
 F_\tau(\rho e_\theta)=\rho
             \frac{p_{r+\tau}(\theta)}{p_r(\theta)}e_\theta
\]
satisfy \eqref{eq:app-map-limits}; their limiting velocity is Lipschitz
because its boundary values match at each fixed vertex ray.
Proposition~\ref{prop:app-volume-derivatives} and the preceding
stress-field argument apply. Differentiating the circle equation
$|x-C(r)|^2=1$ along the moving boundary gives $V\cdot n=C'(r)\cdot n$.
We thus obtain the variation formulas on the full hexagon. No shape
derivative of the auxiliary mixed problem is used.

To include the two endpoints, it remains to prove continuity. The
convex bodies converge in Hausdorff distance and contain a common
disk $B_\rho$, with $\rho>0$. If their Hausdorff distance to a limiting body $D$ is
$\varepsilon$, convexity and $B_\rho\subset D,H_r$ give
\[
 (1+\varepsilon/\rho)^{-1}D\subset H_r
                          \subset(1+\varepsilon/\rho)D.
\]
Domain monotonicity and the scalings $T(tD)=t^4T(D)$ and
$\lambda_1(tD)=t^{-2}\lambda_1(D)$ prove continuity of both functionals.
The same argument gives the Hausdorff continuity used in the angular
maximization and minimization problems.
\end{proof}
\section{First variations for straight caps}
\label{app:straight}
For the cap with axis zero, parameterize the upper tangent segment by
\[
 X_+(s)=r e_\alpha-sJe_\alpha,\qquad
 0<s<\ell=\sqrt{(1-r)^2-r^2},\qquad X_-(s)=\overline{X_+(s)}.
\]
Its outward normal is $e_\alpha$, so $JX_+(s)\cdot e_\alpha=s>0$;
the reflected segment has normal speed $-s$.

We verify that the argument of Appendix~\ref{app:variations} also
applies to straight caps. Using the radial profile
\eqref{eq:straight-radial}, we have, on the cap side of a tangency angle
$\theta_0$,
\[
 p(\theta)=\frac r{\cos(\theta-\theta_0)}
 =r+\frac r2(\theta-\theta_0)^2+O(|\theta-\theta_0|^4).
\]
The value and first angular derivative match the disk. At the apex, its angular derivative has a single jump. We remove this
moving jump by choosing the angular map to follow the apex, as in
\eqref{eq:app-radial-map}.
At an opening cap contact the switching interval has length
$O(\tau)$, the radial remainder is $O(\tau^2)$, and the derivative
remainder is $O(\tau)$. Away from the apex and tangency rays the
active profile is smooth. Thus the maps satisfy every assertion in
\eqref{eq:app-map-limits}, with normal velocity $Jx\cdot n$ on
the moving straight segments and zero on the fixed boundary.

Proposition~\ref{prop:app-volume-derivatives} gives the pulled-back
energy derivatives for both functionals. The stress-divergence
argument in Appendix~\ref{app:variations} applies to these convex
piecewise smooth domains as well, yielding the same boundary formulas.
We therefore obtain the one-sided Hadamard derivatives used in
Theorems~\ref{prop:balance} and~\ref{thm:eigen-balance}, also when two
caps are initially tangent. As in the curved case, bounded and almost-everywhere convergent
derivative quotients suffice; differentiability in the
$W^{1,\infty}$ norm is not required.
\section*{Acknowledgements}

The author acknowledges support through the
ANR project STOIQUES. The author would like to thank the Isaac Newton Institute for Mathematical Sciences, Cambridge, for support and hospitality during the programme Geometric spectral theory and applications, where part of the work on this paper was undertaken. This work was supported by EPSRC grant EP/Z000580/1.

{\bf AI assistance:} GPT 6 Astra and 6.1 Sol were used during the preparation of this manuscript. The author guided the writing process, verified the results and takes full responsibility for them.

\section*{Statements and Declarations}

\noindent\textbf{Funding.}
This work was supported through the Agence nationale de la recherche
(ANR), project STOIQUES, and the Engineering and Physical Sciences
Research Council (EPSRC), grant EP/Z000580/1.

\noindent\textbf{Competing interests.}
The author declares no relevant financial or non-financial competing
interests.

\noindent\textbf{Data and code availability.}
The analytical results use no external datasets. A Python implementation
of the exploratory stability-constant calculation in
Remark~\ref{kr:rem:stability-size}, together with its numerical output,
is provided as supplementary material. These exploratory calculations
are not used in the proofs.

\noindent\textbf{Author contributions.}
Beniamin Bogosel is the sole author and is responsible for all aspects
of the work and the final manuscript.

\end{document}